\documentclass[final,3p,times]{elsarticle}

\usepackage{amssymb,amsmath,amsfonts,mathtools}
\usepackage{amsthm}
	\theoremstyle{remark}
	\newtheorem{rem}{Remark}[section]

	\theoremstyle{plain}

\usepackage{graphicx}
\usepackage{graphbox}
\usepackage[colorinlistoftodos]{todonotes}
\usepackage[colorlinks=true, allcolors=blue]{hyperref}
\usepackage{algorithm}
\usepackage{algpseudocode}
\usepackage{subcaption}
\usepackage{booktabs}
\usepackage{multirow}

\usepackage{pgfplots}
\usepackage{tikz}

\biboptions{sort&compress}
\usepackage{lineno,hyperref}
\usepackage{hypernat}
\modulolinenumbers[5]
\usepackage[]{setspace}

\usepackage{todonotes}

\usepackage{listings}
\usepackage{color} 
\definecolor{mygreen}{RGB}{28,172,0} 
\definecolor{mylilas}{RGB}{170,55,241}

\renewcommand{\vec}[1]{\mathbf{#1}}

\newcommand{\tensorOne}[1]{\boldsymbol{#1}}
\newcommand{\tensorTwo}[1]{\boldsymbol{#1}}
\newcommand{\tensorFour}[1]{\mathbb{#1}}
\newcommand{\linearOperator}[1]{#1}
\newcommand{\funSpace}[1]{#1}
\newcommand{\vecFunSpace}[1]{\boldsymbol{#1}}

\newcommand{\Mat}[1]{#1}
\renewcommand{\Vec}[1]{\mathbf{#1}}

\journal{arXiv}

\begin{document}

\begin{frontmatter}

\title{Enhanced multiscale restriction-smoothed basis (MsRSB) preconditioning with applications to porous media flow and geomechanics\tnoteref{t1}}
\tnotetext[t1]{This document is a collaborative effort.}

\author[St]{Sebastian BM Bosma\corref{cor1}}
\ead{sbosma@stanford.edu}

\author[St]{Sergey Klevtsov\corref{cor1}}
\ead{klevtsov@stanford.edu}

\author[Sintef,NTNU]{Olav M{\o}yner}
\ead{olav.moyner@sintef.no}

\author[LLNL]{Nicola Castelletto}
\ead{castelletto1@llnl.gov}

\cortext[cor1]{Corresponding authors}

\address[St]{Department of Energy Resources Engineering, Stanford University, Stanford, CA, USA}
\address[Sintef]{SINTEF Digital}
\address[NTNU]{Norwegian University of Science and Technology, Department  of Mathematical Sciences}
\address[LLNL]{Atmospheric, Earth, and Energy Division, Lawrence Livermore National Laboratory, Livermore, CA, U.S.A.}

\begin{abstract}
A novel method to enable application of the Multiscale Restricted Smoothed Basis (MsRSB) method to non M-matrices is presented.
The original MsRSB method is enhanced with a filtering strategy enforcing M-matrix properties to enable the robust application of MsRSB as a preconditioner. 
Through applications to porous media flow and linear elastic geomechanics, the method is proven to be effective for scalar and vector problems with multipoint finite volume (FV) and finite element (FE) discretization schemes, respectively. Realistic complex (un)structured two- and three-dimensional test cases are considered to illustrate the method's performance.

\end{abstract}


\end{frontmatter}



\allowdisplaybreaks

\section{Introduction}
Large-scale numerical simulations are often required to understand and predict real world dynamics.
In many applications, the use of high-resolution grids is required to characterize the heterogeneity of the material properties and the geometric complexity of the domains.
Such simulations impose severe computational challenges and motivate the need for efficient solution schemes.
Attractive multilevel strategies to achieve this are multiscale methods \cite{EfeHou09}.
In this paper, we propose a generalization of the multiscale restriction-smoothed basis method (MsRSB) recently put forward in \cite{MsRSB_Moyner2016}, and investigate its use as an effective preconditioner for multipoint flux approximation finite volume (FV) and finite element (FE) discretizations of second-order elliptic problems. Specifically we focus on applications to porous media flow and linear elastic geomechanics.

The original idea underlying multiscale discretization methods for heterogeneous second-order elliptic problems can be traced back four decades \cite{StrFed79,BabOsb83}.
In essence, these methods aim at constructing accurate coarse-scale problems that preserve information of fine scale heterogeneity and can be solved at low computational cost.
This is accomplished by numerically computing multiscale basis functions, which are local solutions of the original problem, that are used to both: (i) construct the coarse-scale problem, and (ii) interpolate the coarse-scale solution back to the fine-scale. 
Various methods to obtain these basis functions have been developed, for example generalized finite-element (GFE) methods \cite{BabCalOsb94}, multiscale finite-element (MsFE) methods \cite{HouWu97}, numerical-subgrid upscaling \cite{Arb02}, multiscale mixed finite-element (MsMFE) methods \cite{CheHou03}, multiscale finite-volume (MsFV) methods \cite{MSFV_Jenny2003},
multiscale mortar mixed finite-element (MsMMFE) methods \cite{Arb_etal07}, multilevel multiscale mimetic (M\textsuperscript{3}) methods \cite{LipMouSvy08}, multiscale mixed/mimetic finite-element (MsMFEM) \cite{AarKrogLie08} and generalized multiscale finite element (GMsFE) \cite{EfeGalHou13} methods, to name a few. 
In the geoscience community, multiscale methods have been extensively applied both as single-pass \cite{MSFV_Jenny2003} and iterative schemes \cite{NorBjo08,Haj_etal08} to resolve some of the limitations of existing upscaling methods.
They have established a solid framework for simulating complex subsurface flow processes, e.g. \cite{JuaTch08,Nor09,Hel_etal10,ZhoTch12,WanHajTch14,Koz_etal16,Cus_etal15,ChuEfeLee15,ParEdw15,TenKobHaj16,Lie_etal16,LieMoyNat17,Cus_etal18}.
Multiscale methods for linear elastic problems have focused primarily on the derivation of accurate coarse space basis functions which are robust with respect to material property heterogeneities and enable scalable performance \cite{BucIliAnd13,BucIliAnd14,Spi_etal14,MultiscaleFEM_Castelleto2017,ChuLee19}.
Applications to the poroelasticity equations include \cite{ZhaFuWu09,BroVas16a,BroVas16b,DanGanWhe18,Akk_etal18,SokBasHaj19,Cas_etal19}.

The MsRSB method was proposed in the context of FV simulation for fluid flow in highly heterogeneous porous media \cite{MsRSB_Moyner2016}.
Based on a two-grid approach, the MsRSB method constructs multiscale basis functions through restricted smoothing on the fine-scale matrix.
In more detail, the basis functions, which are consistent with the local differential operators, are constructed with a cheap relaxation scheme, i.e. a weighted Jacobi iteration, similar to approaches used in smoothed aggregation multigrid methods \cite{VanManBre96,VanManBre01,Bre_etal05}.
An important advantage of MsRSB is that smoothing by relaxation provides a great deal of flexibility in handling unstructured grids, an essential requirement, for example, in applications involving complex geological structures.
MsRSB has been widely proven and implemented in open source and commercial simulators using a linear two-point flux approximation (TPFA) \cite{Lie_etal16}.

Because of the two-point structure, the linear TPFA scheme is monotone \cite{Dro14}, i.e. it preserves the positivity of the differential solution \cite{BerPle94}, and leads to an M-matrix with a small stencil.
This is the reason why linear TPFA is the scheme of choice in most engineering software.
Unfortunately, the consistency of TPFA is not guaranteed for arbitrary grids and anisotropic permeability distributions, potentially leading to inaccurate results \cite{MRST}. Therefore, other FV methods such as multipoint flux approximation (MPFA) and/or nonlinear schemes \cite{Dro14,TerMalTch17} must be considered to achieve consistent fluxes.
To date, few works have investigated MsRSB applied to MPFA or other consistent discretizations \cite{souza20}. 
Moreover, to the authors' knowledge, the issues associated to non M-matrices have not been addressed in the literature.
Hence, in this paper, we focus on enhancing MsRSB to enable the solution of second-order elliptic problems using discretization methods that do not result in an M-matrix, thus allowing general application to consistent discretizations.
Based on the MPFA-O method \cite{MPFA_Aavatsmark}, we show that the MsRSB basis construction as presented in \cite{MsRSB_Moyner2016} can fail due to divergent iterations for an anisotropic diffusion problem.
We propose a variant of the original MsRSB approach that restores the desired behavior by enforcing M-matrix properties based on a filtering strategy.
We develop the new method focusing on FV discretizations for porous media single-phase flow, and extend its use to vector elliptic problems by targeting FE-based simulation of linear elastic geomechanics.

The paper is structured as follows.
First, the original multiscale restriction-smoothed basis method is briefly reviewed in \ref{sec:MsRSB}.
Second, MsRSB for an MPFA flow discretization is analyzed and the novel approach is proposed in Section \ref{sec:MPFA}.
Next, the proposed method is extended to geomechanics in Section \ref{sec:geomechanics}.
Challenging two- and three-dimensional experiments are presented to demonstrate properties, robustness and scalability of the method throughout Section \ref{sec:MPFA} and \ref{sec:geomechanics}, including comparisons to existing methods and published results.
Finally the report is concluded and future work specified.

\section{The Multiscale Restriction-Smoothed Basis method (MsRSB)}  \label{sec:MsRSB}
We propose a two-level preconditioning framework based on MsRSB for accelerating iterative Krylov methods to solve linear systems of the form: 
\begin{linenomath}
\begin{equation}
  \Mat{A} \vec{u}  = \vec{f},
  \label{eq:linsys_general}
\end{equation}
\end{linenomath}
where the coefficient matrix $\Mat{A} \in \mathbb{R}^{n \times n}$ arises from a finite volume (FV) or finite element (FE) discretization of a scalar or vector second-order elliptic problem. Furthermore,  $\vec{u} = \{ u_i \}_{i=1}^{n} \in \mathbb{R}^{n}$ is the solution vector containing the unknown degrees of freedom, and $\vec{f} = \{ f_j \}_{j=1}^{n} \in \mathbb{R}^{n}$ is the discrete forcing term.
In this work we develop the method and illustrate its performance focusing on two simple but representative models routinely employed in practical simulation of subsurface processes: (i) the incompressible single-phase flow equation, and (ii) the linear elastostatic equations.
For the flow problem we will concentrate on FV fine-scale discretizations while for the elastostatics problem we will consider the FE method.
A review of governing equations and the derivation of the matrix form in \eqref{eq:linsys_general} are provided for both models in \ref{app:models_strong_form} and \ref{app:model_problems_discretization}, respectively.

The essence of any multiscale formulation and solution algorithm lies in the construction of a representative coarse-scale problem capable of capturing fine-scale features of a high-resolution model.
The connection among scales is accounted for by computing basis functions, i.e. localized fine-scale solutions, which are used to construct a coarse-scale (upscaled) problem and reconstruct a fine-scale (downscaled) solution from the coarse solution.
The reconstruction stage can be represented through the \textit{prolongation operator} $\linearOperator{P}$, a sparse linear operator that stores the basis function associated to each coarse degree of freedom in the corresponding column such that

\begin{align}
  \linearOperator{P} : \mathbb{R}^{n_c} \rightarrow \mathbb{R}^{n_f}, \vec{u}_c \mapsto \vec{u}_f = \linearOperator{P}\vec{u}_c.
  \label{eq:prolongation}
\end{align}

\noindent
Here and in the following, subscripts $f$ and $c$ indicate quantities associated with the fine and the coarse problem.
In particular, $n_c$ and $n_f$ denote the number of coarse- and fine-scale degrees of freedom, respectively.

Assuming $\linearOperator{P}$, to be specified hereafter, is available, the definition of the coarse problem proceeds as follows.
First the fine-scale solution in \eqref{eq:linsys_general} is replaced with the approximation $\vec{u} = \vec{u}_f \approx \linearOperator{P} \vec{u}_c$.
Second the resulting residual vector, namely $\vec{r} = (\vec{f} - \Mat{A} \linearOperator{P} \vec{u}_c)$, is orthogonalized against $n_c$ vectors in $\mathbb{R}^{n_c}$ that form the rows of the operator $\linearOperator{R}$.
Hence, $\vec{u}_c$ is the solution to the linear system with $n_c$ equations

\begin{align}
  \Mat{A}_c \Vec{u}_c &= \Vec{f}_c &
  &\text{with} &  
  \Mat{A}_c &= \linearOperator{R} \Mat{A} \linearOperator{P}, &
  \vec{f}_c &= \linearOperator{R} \vec{f}.
  \label{eq:lin_system_coarse}
\end{align}

\noindent
We refer to $\linearOperator{R}$ as the \textit{restriction operator}, a linear operator mapping vectors from the fine- to the coarse-scale

\begin{align}
  \linearOperator{R} : \mathbb{R}^{n_f} \rightarrow \mathbb{R}^{n_c}, \vec{u}_f \mapsto \vec{u}_c = \linearOperator{R}\vec{u}_f.
  \label{eq:restriction}
\end{align}

\noindent
Different options may be considered for $\linearOperator{R}$.
If the fine-scale matrix $\Mat{A}$ is symmetric positive definite (SPD), a Galerkin orthogonalization, i.e. $\linearOperator{R} = \linearOperator{P}^T$, is typically the strategy of choice since it provides a coarse scale operator $\Mat{A}_c$ that is still SPD.
Alternatively, a Petrov-Galerkin approach is often used.
For example, in the MSFV method \cite{MSFV_Jenny2003}, which is designed for diffusion problems, $\linearOperator{R}$ is constructed such that discrete mass conservation also holds for the coarse-scale problem.

\subsection{MsRSB for TPFA finite volume schemes} \label{sec:FvTpfaMsRSB}
\begin{figure} [htbp]
\centering
\begin{subfigure}[t]{0.3\textwidth}
\centering
\includegraphics[width=\linewidth]{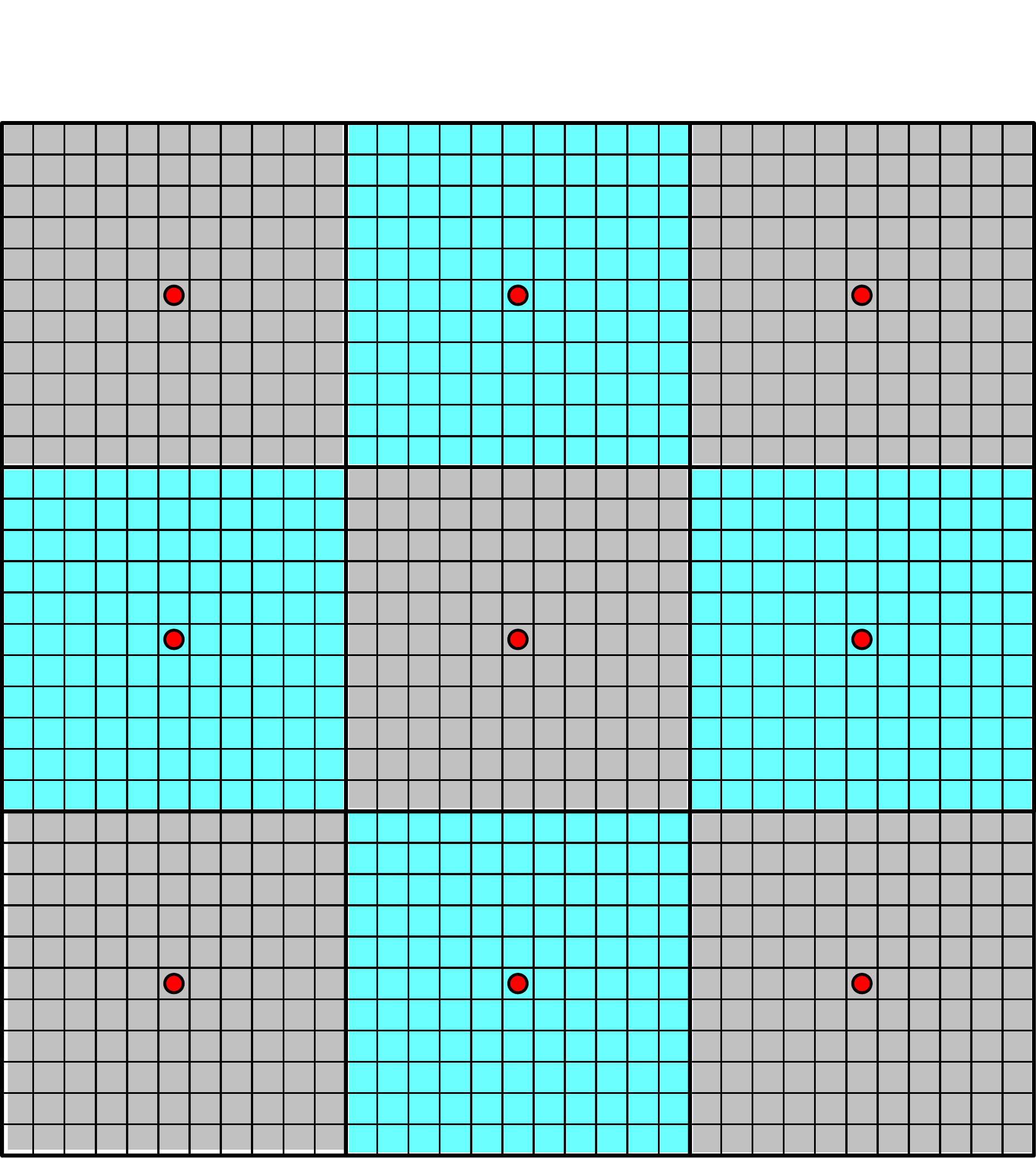}\hfill 
\caption{\label{fig:OlavGrid_new_a}}
\end{subfigure}
\hfill
\begin{subfigure}[t]{0.3\textwidth}
\centering
\includegraphics[width=\linewidth]{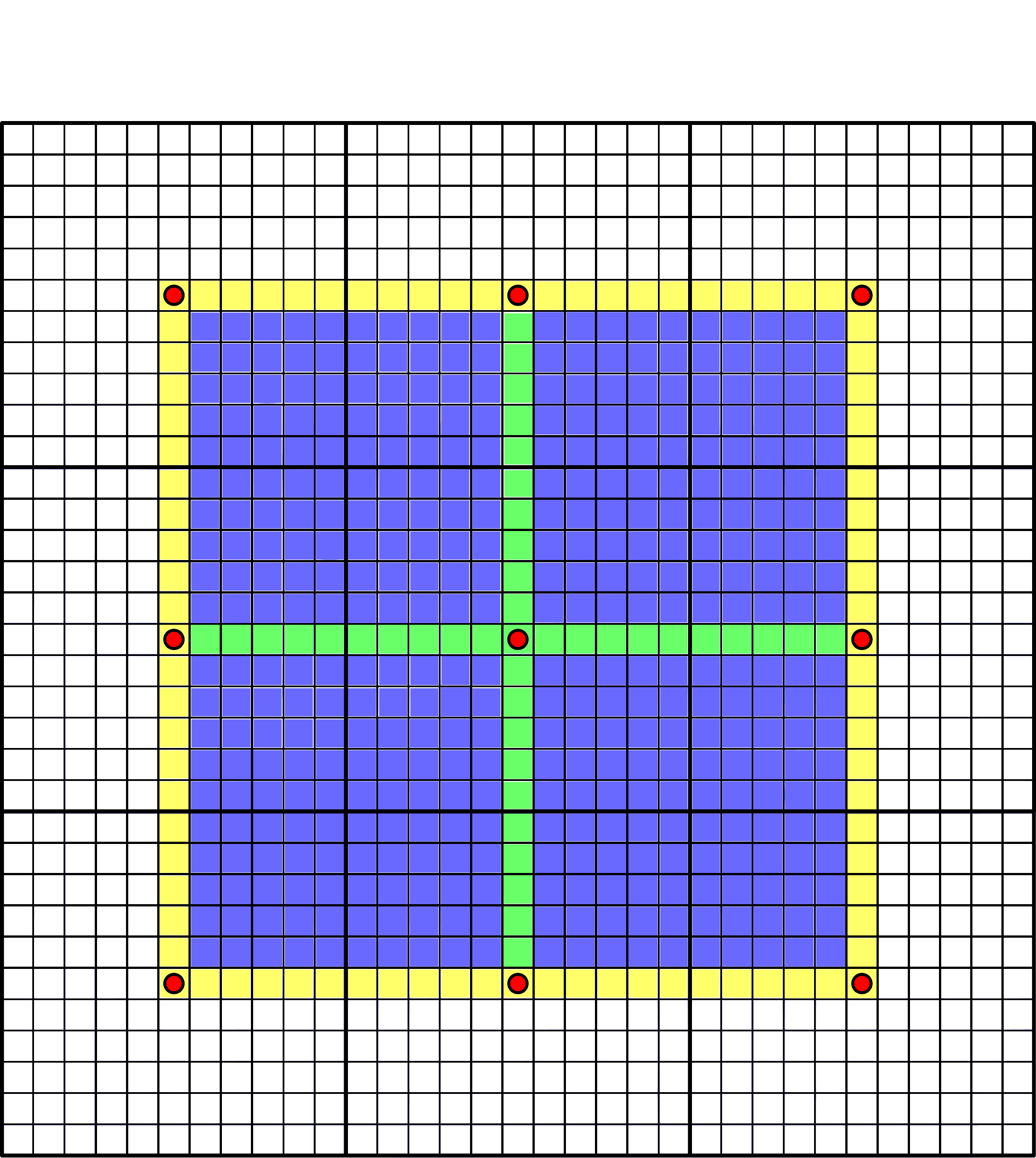}\hfill 
\caption{\label{fig:OlavGrid_new_b}}
\end{subfigure}
\hfill
\begin{subfigure}[t]{0.3\textwidth}
\centering
\includegraphics[width=\linewidth]{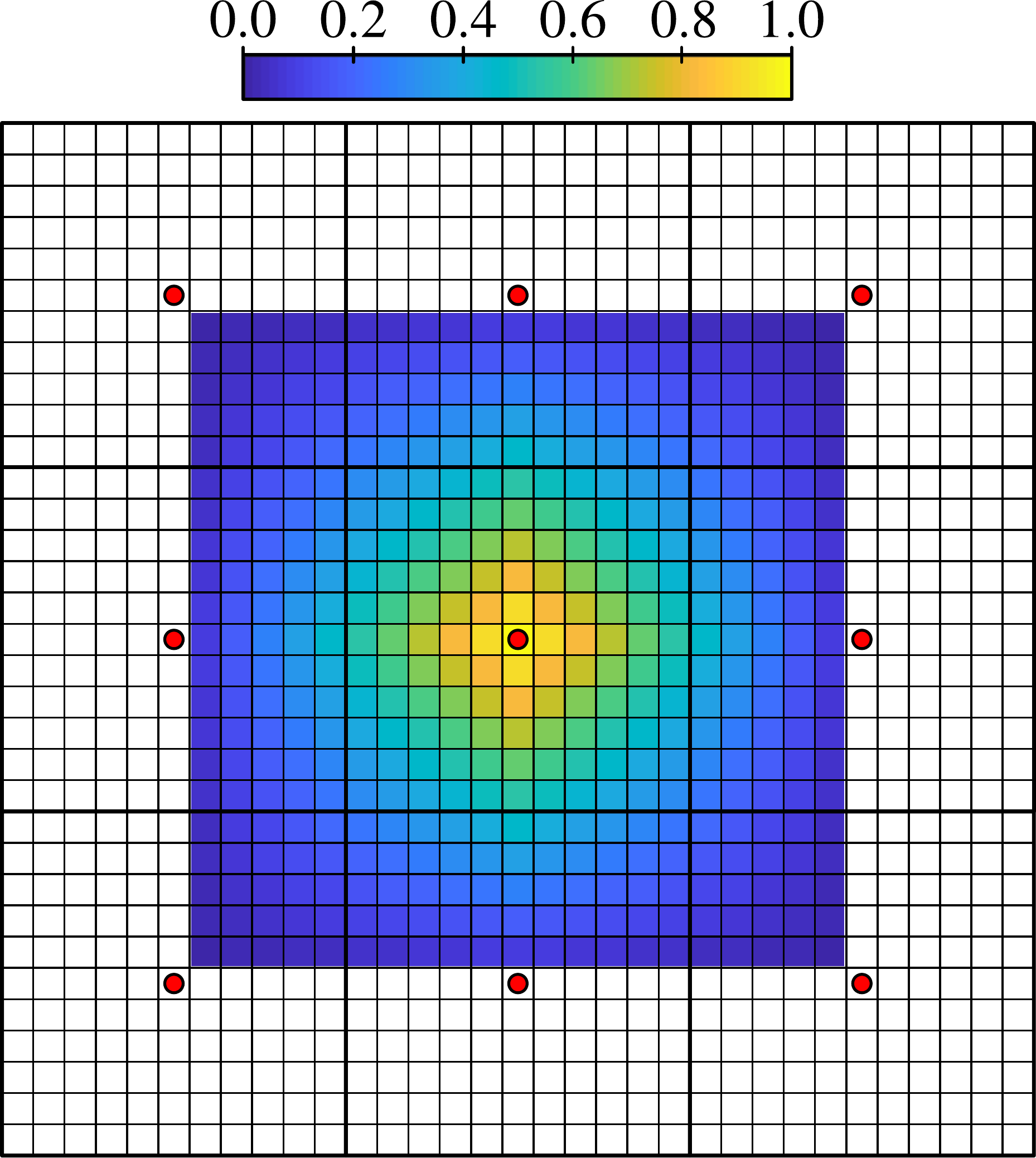}\hfill 
\caption{\label{fig:OlavGrid_new_c}}
\end{subfigure}
\caption{\label{fig:OlavGrid_new} Coarse grid features for the internal block of a regular 3$\times$3 partition: (a) primal grid and coarse nodes (\protect\tikz \protect\draw[black,fill=red] (0,0) circle (.5ex););  (b) support boundary (\protect\tikz \protect\draw[black,fill=yellow] (0,0) rectangle (1ex, 1ex);), and edge (\protect\tikz \protect\draw[black,fill=green] (0,0) rectangle (1ex, 1ex);) cells; and (c) basis function.}
\end{figure}

A crucial component of multiscale methods is the efficient and accurate construction of the prolongation operator, that is the computation of the basis functions.
Originally proposed for the cell-centered finite volume solution to the diffusion equation for flow through porous media \cite{MsRSB_Moyner2016}, the MsRSB method computes the basis functions iteratively with restricted smoothing.
To describe this process, some terminology is first defined.
A primary coarse grid is defined as a partitioning of the fine grid.
In each coarse cell, a coarse node is chosen as the representative fine cell for that coarse cell.
Furthermore, support boundary cells are defined as the cells connecting neighboring coarse nodes.
Support edge cells are defined as the cells connecting a coarse node to its neighbors.
Note that edge cells for one coarse node will be boundary cells for other coarse nodes.
Fig. \ref{fig:OlavGrid_new} shows an example of the coarse grid structure for a flow problem.

The initial guess for the basis function associated to the $j$th coarse node consists of the characteristic function of the primary $j$th coarse cell, hence it is equal to 1 for fine scale cells belonging to the $j$th coarse cell and 0 elsewhere.
\newcommand{\basismat}{\ensuremath{\Mat{G}}}
Let $\Mat{L}$ and $\Mat{U}$ denote the strictly lower-triangular and upper-triangular part of the matrix \basismat{} which we want to compute basis functions for. We assume that \basismat{} has the same dimensions as the fine-scale discretization matrix $\Mat{A}$, but otherwise leave the relationship ambiguous for the time being.
Let $\hat{\Mat{D}}$ be the diagonal matrix such that the entries in each row of matrix $\hat{\basismat} = (\Mat{L} + \hat{\Mat{D}} + \Mat{U})$ sum to zero, i.e.
\begin{linenomath}
\begin{align}
  [\hat{\Mat{D}}]_{ii} = - \sum_{j=1,j \ne i}^{n} [\basismat]_{ij}, \qquad \forall i \in \{1, 2 \ldots, n \}.
  \label{eq:Dhat_zero_rowsum}
\end{align}
\end{linenomath}
Each MsRSB basis function is then iteratively smoothed by applying the following relaxation scheme

\begin{linenomath}
\begin{align}
  [\linearOperator{P}]_{*j}^{k+1} = [\linearOperator{P}]_{*j}^k - \omega \hat{\Mat{D}}^{-1} \hat{\basismat} [\linearOperator{P}]_{*j}^k,   
  \label{eq:smoother}
\end{align}
\end{linenomath}
where $[\linearOperator{P}]_{*j}$ denotes the $j$th column of $\linearOperator{P}$, i.e. the basis function corresponding to coarse node $j$, and $k$ is the iteration count.
Furthermore, $\omega$ is a relaxation factor which in this work is set to 2/3, i.e. the value warranting the optimal smoothing factor of the weighted Jacobi iteration for the homogeneous Poisson's equation  \cite{Saa03}.
If the smoothing update extends the basis function outside of its support region, the update is adjusted to enforce that the basis support is enclosed by the boundary cells.
This is done by adding overflowing basis function values to the neighboring basis functions.
We refer the reader to \cite{MsRSB_Moyner2016} for additional details.
An important feature of the originally described smoothing process is that it theoretically guarantees the conservation of the initial partition of unity at each iteration. However, round-off errors in a numerical implementation can cause a violation of this condition. An efficient and robust MsRSB implementation, which guarantees the numerical partition of unity, is described in algorithm 2 of \cite{JohannessenThesis}. In short, the updated basis function values are rescaled rather than the updates.
Section 3.7 of the thesis elaborates on the issues of round off error.
The full MsRSB basis function construction process is also incorporated in algorithm \ref{alg:EnhancedMsRSB}.


\section{MsRSB for non-M matrices: multipoint FV schemes} \label{sec:MPFA}
In order to demonstrate the challenges of applying a multiscale method to discretized systems which do not result in M-matrices, we will consider a specific fine-scale discretization of the incompressible single-phase pressure equation 

\begin{equation}
    - \text{div}(\tensorTwo{\Lambda} \cdot \text{grad} \; p) = q,
    \label{eq:pressure}
\end{equation}

\noindent
where $p$ is the scalar pressure field, $q$ a source term distributed in the domain and $\tensorTwo{\Lambda}$ a positive-definite tensor describing the diffusion properties of the medium. Generally, $\tensorTwo{\Lambda}$ is characterized by the permeability tensor $\tensorTwo{\kappa}$ and fluid viscosity $\mu$.
A detailed description of the governing equations is provided in \ref{app:model_flow}.
The standard approach for discretizing \eqref{eq:pressure} is a two-point flux approximation (TPFA) scheme, which is only consistent when the principal axes of the permeability tensor are aligned with the grid \cite{EdwRog98}---i.e., for so-called $\tensorTwo{\kappa}$-orthogonal grids. 
One possible choice to solve pressure on rough grids with non-diagonal permeability tensors is the MPFA-O method \cite{MPFA_Aavatsmark}, which may not result in M-matrices for grids of interest.
In the MPFA-O method, transmissibilities are computed by enforcing continuity of fluxes over each half-face for local reconstructions of linear flow.
These fluxes rely on all the cell pressure unknowns surrounding a grid vertex and guarantee by construction that the set of equations, related to the half faces connected to a grid vertex, result in a solvable system.
As it is not the focus of this paper, the authors refer to the cited paper for more details on the well-established method.
Furthermore, we note that the implementation of MsRSB for MPFA was done using the Matlab Reservoir Simulation Toolbox \cite{MRST}. 
In the following, we address two points: i) we will consider scalar systems discretized with MPFA-O as a proxy for the inherent difficulties in applying multiscale methods to non-M-matrices and ii) demonstrate a practical approach for implementing MsRSB or similar methods to MPFA-type discretizations for flow.
%

\subsection{Extension to non M-matrices}
\label{sec:fix_non_M-matrix}
To avoid complications related to non M-matrices, a method is devised to alter, with minimal intrusion, the original fine-scale matrix $\Mat{A}$
such that it satisfies M-matrix properties.
This approximation of the linear system is justified because we aim to find an approximate solution using multiscale methods, where the basis functions should account for the M-matrix like part of the system matrix.
This is conceptually similar to e.g. using incompressible or steady-state basis functions for flow when the problem under consideration is nonlinear due to compressibility \cite{Compressible_MS_Zhou,MsRSB_compressible_Moyner2016}.
Moreover, this work's primary objective is to employ the proposed method as a preconditioner for GMRES and other iterative Krylov solvers, where an inexpensive local solver will target any local errors in the approximation. 
In the following, we assume that the discretization matrix has by convention a positive diagonal and primarily negative off-diagonal entries.

\begin{figure} [htbp]
\begin{subfigure}[t]{0.22\textwidth}
  \centerline{\includegraphics[width=\linewidth]{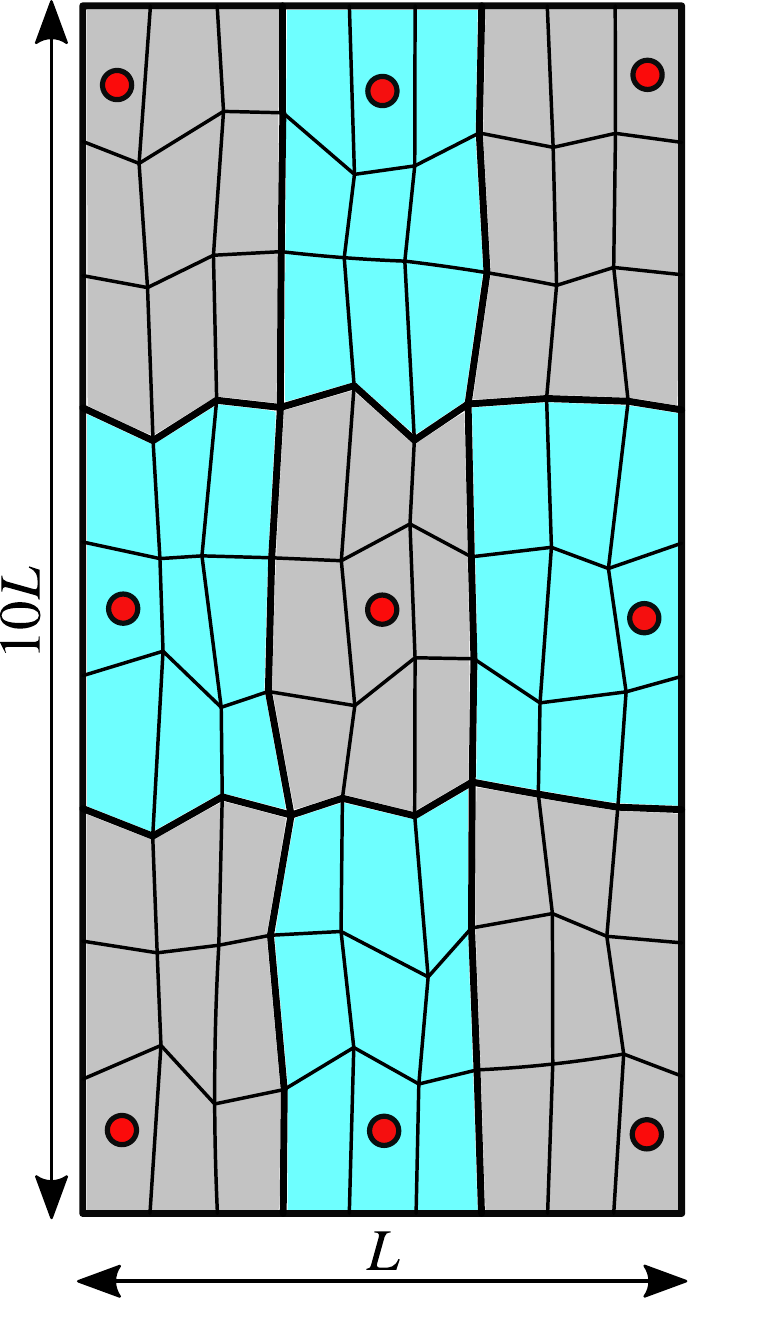}}
  \caption{\label{fig:CoarseGridMPFA}}
\end{subfigure}
\hfill
\begin{subfigure}[t]{0.22\textwidth}
  \centerline{\includegraphics[width=\linewidth]{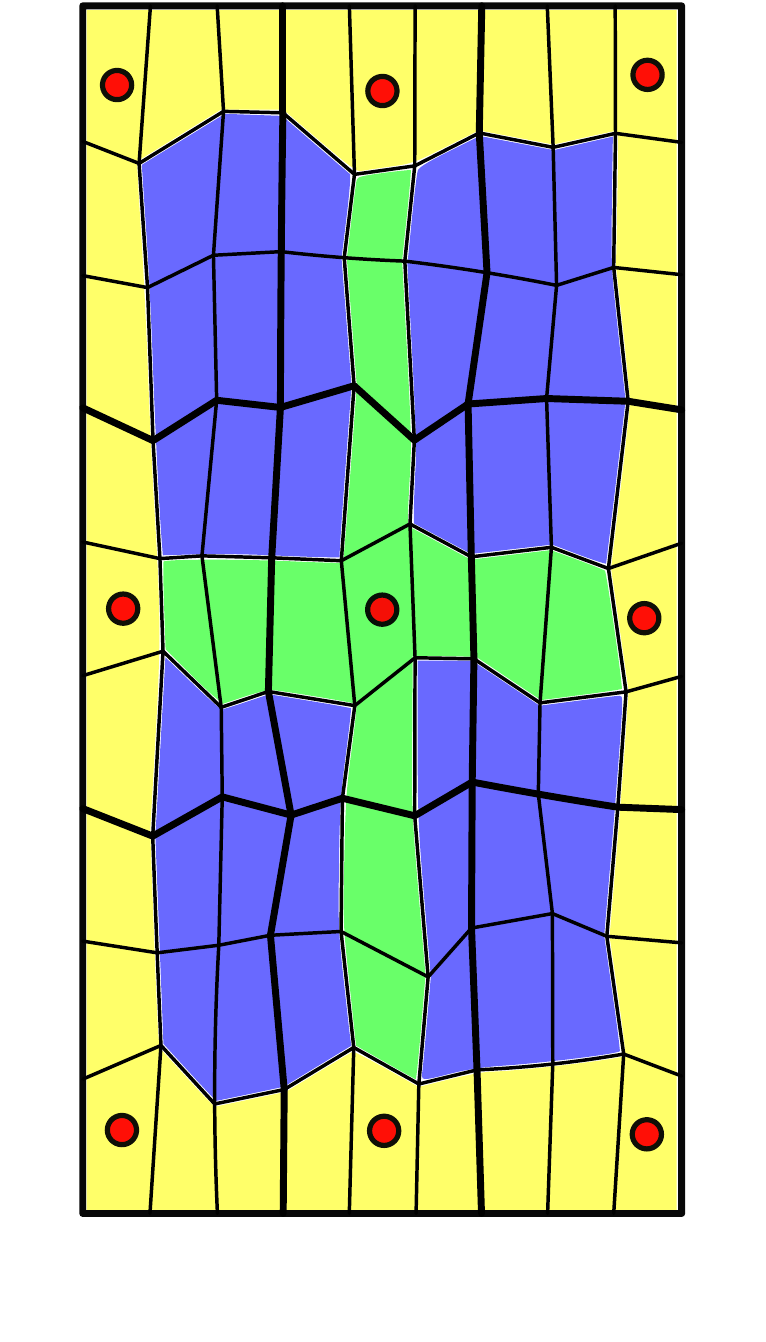}}
  \caption{\label{fig:MPFAb}}
\end{subfigure}
\hfill
\begin{subfigure}[t]{0.22\textwidth}
  \centerline{\includegraphics[width=\linewidth]{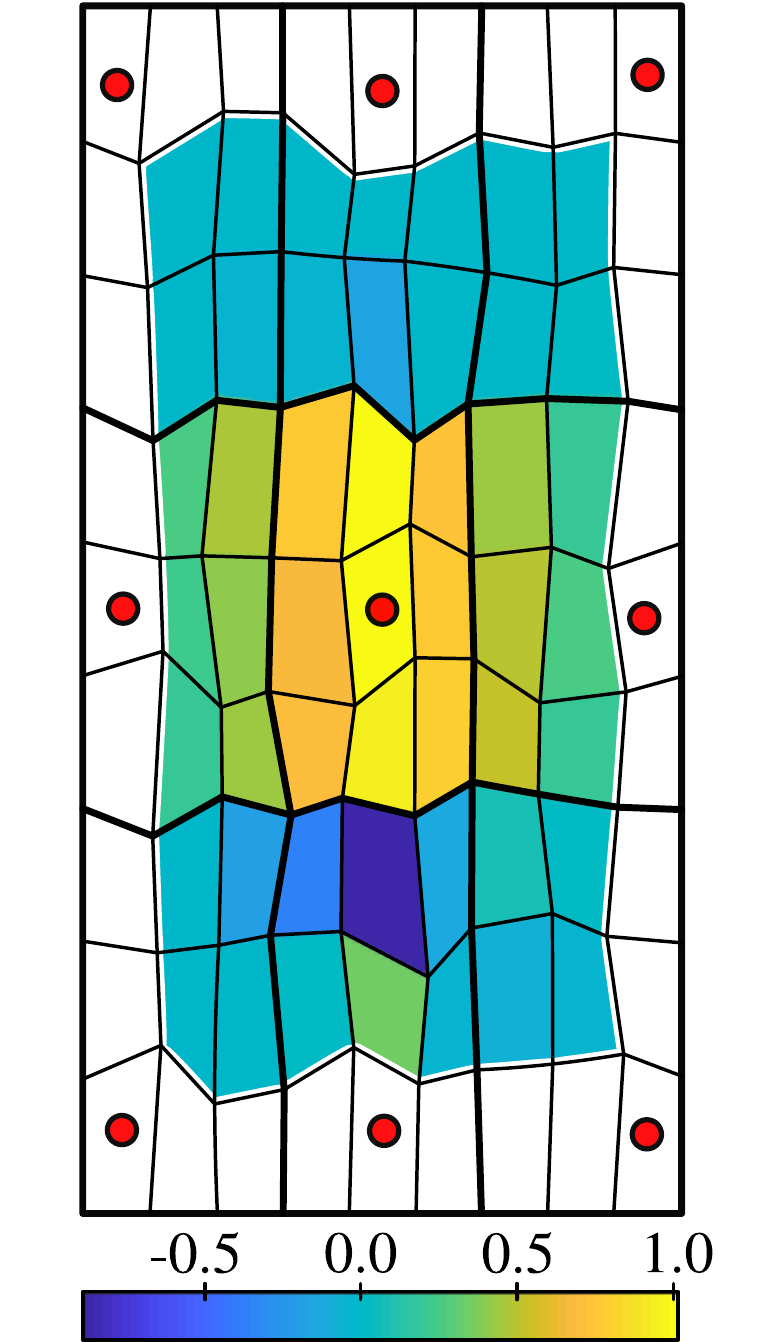}}
  \caption{\label{fig:MPFAbforig}}
\end{subfigure}
\hfill
\begin{subfigure}[t]{0.22\textwidth}
  \centerline{\includegraphics[width=\linewidth]{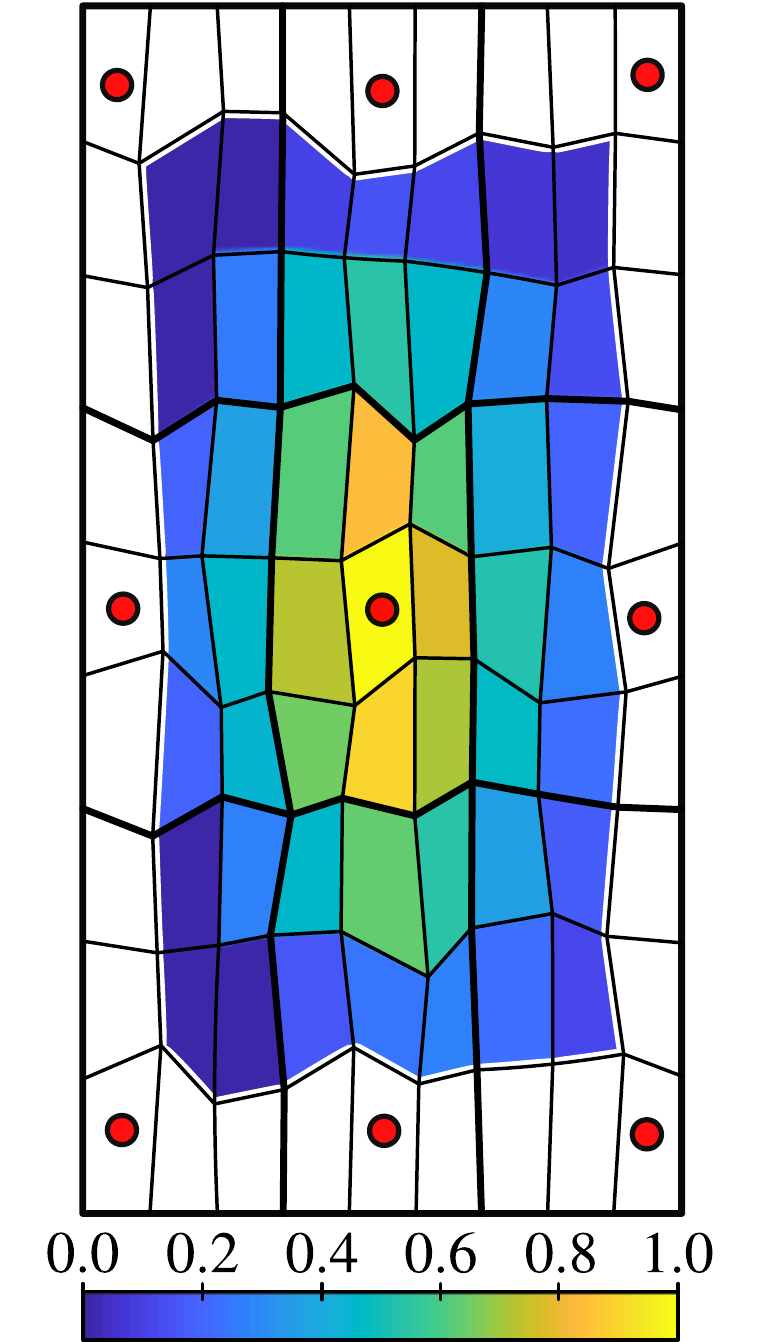}}
  \caption{\label{fig:MPFAbfalter}}
\end{subfigure}
\caption{\label{fig:MPFAbfs} Basis function for the internal block of a regular 3$\times$3 partition: (a) primal grid and coarse nodes (\protect\tikz \protect\draw[black,fill=red] (0,0) circle (.5ex);); (b) support boundary (\protect\tikz \protect\draw[black,fill=yellow] (0,0) rectangle (1ex, 1ex);), edge (\protect\tikz \protect\draw[black,fill=green] (0,0) rectangle (1ex, 1ex);) and internal cells (\protect\tikz \protect\draw[black,fill=blue] (0,0) rectangle (1ex, 1ex);); (c) MsRSB basis function using the original fine-scale system; and (d) MsRSB basis function using the filtered fine-scale system. The convergence criteria is $e_{it} < 10^{-12}$, where $e_{it}$ is the maximum value in the smoother update on the internal cells (see algorithm \ref{alg:EnhancedMsRSB}). However as the case with the original linear system diverges, the basis function obtained after 11 iterations is plotted in (c).}
\end{figure}

To enforce M-matrix properties, it is sufficient to filter out all positive off-diagonal entries and construct a modified system matrix $\tilde{\Mat{A}}$
\begin{linenomath}
\begin{align}
    [\tilde{\Mat{A}}]_{ij} = min([\Mat{A}]_{ij},0), \qquad \forall (i,j) \in \{1, 2, \ldots, n \} \times \{1, 2, \ldots, n \},
 \label{eq:negoffdiag}
\end{align}
\end{linenomath}
from which we can easily compute the basis functions by setting $\basismat{} = \tilde{\Mat{A}}$. As the resulting matrix $\hat{\basismat{}}$ from eq. \eqref{eq:smoother} now has zero row sum with only non-positive off-diagonal entries, the enhanced MsRSB method guarantees the robust generation of partition-of-unity basis functions with entries in $[0, 1]$.
Furthermore, note that these changes do not affect a problem described by an M-matrix and as such the method can be implemented generally.
To verify that the method has the desired effect, a simple test case is devised. Starting from an equidistant 2-D Cartesian grid, all grid vertices are perturbed randomly in both $x$- and $y$-direction. Additionally the grid is stretched by a factor 10 in the $y$-direction. The resulting test case primal grid is shown in Fig. \ref{fig:CoarseGridMPFA}. The fine grid has 9x9 cells where we coarsen by a ratio of 3 in each direction. Furthermore all flow properties are homogeneous. Note that MPFA-O for a $\tensorTwo{\kappa}$-orthogonal Cartesian grid is equivalent to the TPFA method. Therefore the test case is designed to vary from an orthogonal grid substantially.

The basis functions corresponding to the central coarse node are plotted in Figure \ref{fig:MPFAbfs}.
Figure \ref{fig:MPFAbforig} displays the basis function obtained after 11 iterations when using the original linear system.
It is evident from the plot that the prolongator is diverging. Figure \ref{fig:MPFAbfalter} presents the converged basis function when using the altered linear system. In this case, the desired monotone basis functions are obtained. Note that due to the cancellation of certain connections in the matrix, the basis function does not spread over the full domain. This effect is especially present in the corners of the dual-region. As will be shown in numerical examples in Section \ref{sec:geomechanics}, although this has no impact on the final solution, it is favorable to avoid small coarsening ratios to ensure good performance. The proposed multiscale method is summarized in algorithm \ref{alg:EnhancedMsRSB}. Furthermore, the multiscale preconditioning strategy employed in the numerical results is described in algorithms \ref{alg:setupPreconditioner}-\ref{alg:applyPreconditionerTwoLevel}.


\begin{algorithm}
\caption{Enhanced MsRSB method to construct prolongation operator.} \label{alg:EnhancedMsRSB}
\begin{algorithmic}[1]
\Function{EnhancedMsRSB}{$\Mat{A}$,$\linearOperator{P}$}
%
%
\State $\hat{\Mat{A}} = \min(\Mat{A}, 0)$ \Comment{Eq. \ref{eq:negoffdiag}}
\State $\hat{\Mat{A}} \gets \hat{\Mat{A}} - \text{diag}(\text{rowsum}(\hat{\Mat{A}}))$ \Comment{Eq. \ref{eq:Dhat_zero_rowsum}}

\State $\hat{\Mat{D}} = \text{diag}(\hat{\Mat{A}})$

\State k = 0, $e_{it} = \infty$

\While{$e_{it} > \text{tol}$}

\State $\delta\linearOperator{P} = - \omega \hat{\Mat{D}}^{-1} \hat{\Mat{A}} \linearOperator{P}$  \Comment{Eq. \ref{eq:smoother}}
\State Modify $\delta\linearOperator{P}$ to avoid stencil growth outside of support region \Comment{See section 3.3 step 2 of \cite{MsRSB_Moyner2016}}
\State $\linearOperator{P} \gets \linearOperator{P} + \delta\linearOperator{P} $
\For {$i = 1:n_{f}$} 
\State $\linearOperator{[P]}_{i*} \gets \linearOperator{[P]}_{i*} / \sum_j(\linearOperator{[P]}_{ij}) $ \Comment{Rescale to guarantee numerical partition of unity}
\EndFor
\If{($k \bmod n_{it}$) = 0} \Comment{Check convergence every $n_{it}$ iterations}
\State $e_{it} = \max\limits_{i,j} (\text{abs}( [\delta \linearOperator{P}]_{ij})), \qquad  i \not \in$ support edges
\EndIf
\State $k \leftarrow k+1$
\EndWhile
\State \textbf{return} $\linearOperator{P}$ \Comment{The converged prolongation operator}
\EndFunction
\end{algorithmic}
\end{algorithm}


\begin{algorithm}
\caption{Setup of enhanced MsRSB preconditioner.} \label{alg:setupPreconditioner}
\begin{algorithmic}[1]
\Procedure{SetupMsRSBpreconditioner}{$\Mat{A}$}
\If {Not initialized}
\State Initialize $\linearOperator{P}$
\EndIf
\State $\linearOperator{P}$ = \Call{EnhancedMsRSB}{$\Mat{A}$,$\linearOperator{P}$} 
\If {Petrov-Galerkin restriction}
  \State ${\linearOperator{R}} \gets$ Construct restriction operator \Comment{Compute restriction operator}
\Else
  \State ${\linearOperator{R}} = {\linearOperator{P}}^T$
\EndIf
\State $\Mat{A}_c = \linearOperator{R} \Mat{A} \linearOperator{P}$
\Comment{Compute coarse-scale system matrix}
\State $\Mat{M}_c^{-1} \approx A_c^{-1}$ \Comment{Set up coarse system solver} 
\EndProcedure
\end{algorithmic}
\end{algorithm}


\begin{algorithm}
\caption{Application of the enhanced MsRSB} \label{alg:applyPreconditioner}
%
%
%
\begin{algorithmic}[1]
\Function{ApplyMsRSBpreconditioner}{$\Vec{v}$}
\State $\Vec{v_c} = \linearOperator{R} \Vec{v}$\
\Comment{Restrict residual}
\State $\Vec{w}_c = M_c^{-1}\Vec{v}_c$
\Comment{Solve coarse problem}
\State $\Vec{w} = \linearOperator{P}\Vec{w}_c$
\Comment{Interpolate the solution}
\State \Return $\Vec{w}$
\EndFunction
\end{algorithmic}
\end{algorithm}
\begin{algorithm}
\caption{Application of the two-level preconditioner with pre- and post-smoothing.} \label{alg:applyPreconditionerTwoLevel}
\begin{algorithmic}[1]
\Function{ApplyMsRSBTwoLevelPreconditioner}{$\Mat{A}$,$\Vec{v}$}
\State $\Vec{z} \gets \vec{z} + M_{\text{pre.}}^{-1} (\Vec{v} - \Mat{A} \Vec{z})$ \Comment{Pre-smoothing relaxation starting from $\Vec{z} = \Vec{0}$ }
\State $\Vec{w}$ = \Call{ApplyMsRSBpreconditioner}{$\Vec{v} - \Mat{A} \Vec{z}$}
\Comment{Compute correction}
\State $\Vec{z} \gets \vec{z} + \Vec{w}$ \Comment{Apply correction}
\State $\Vec{z} \gets \vec{z} + M_{\text{post.}}^{-1} (\Vec{v} - \Mat{A} \Vec{z})$ \Comment{Post-smoothing relaxation}
\State \textbf{return} $\Vec{z}$
\EndFunction
\end{algorithmic}
%
%
%
\end{algorithm}

\subsection{MsRSB for an MPFA-O discretization: 2D test case}
\begin{figure} [htbp]
\centering
 \begin{subfigure}{0.45\textwidth}
 \centering
 \hfill
 \includegraphics[width=\linewidth]{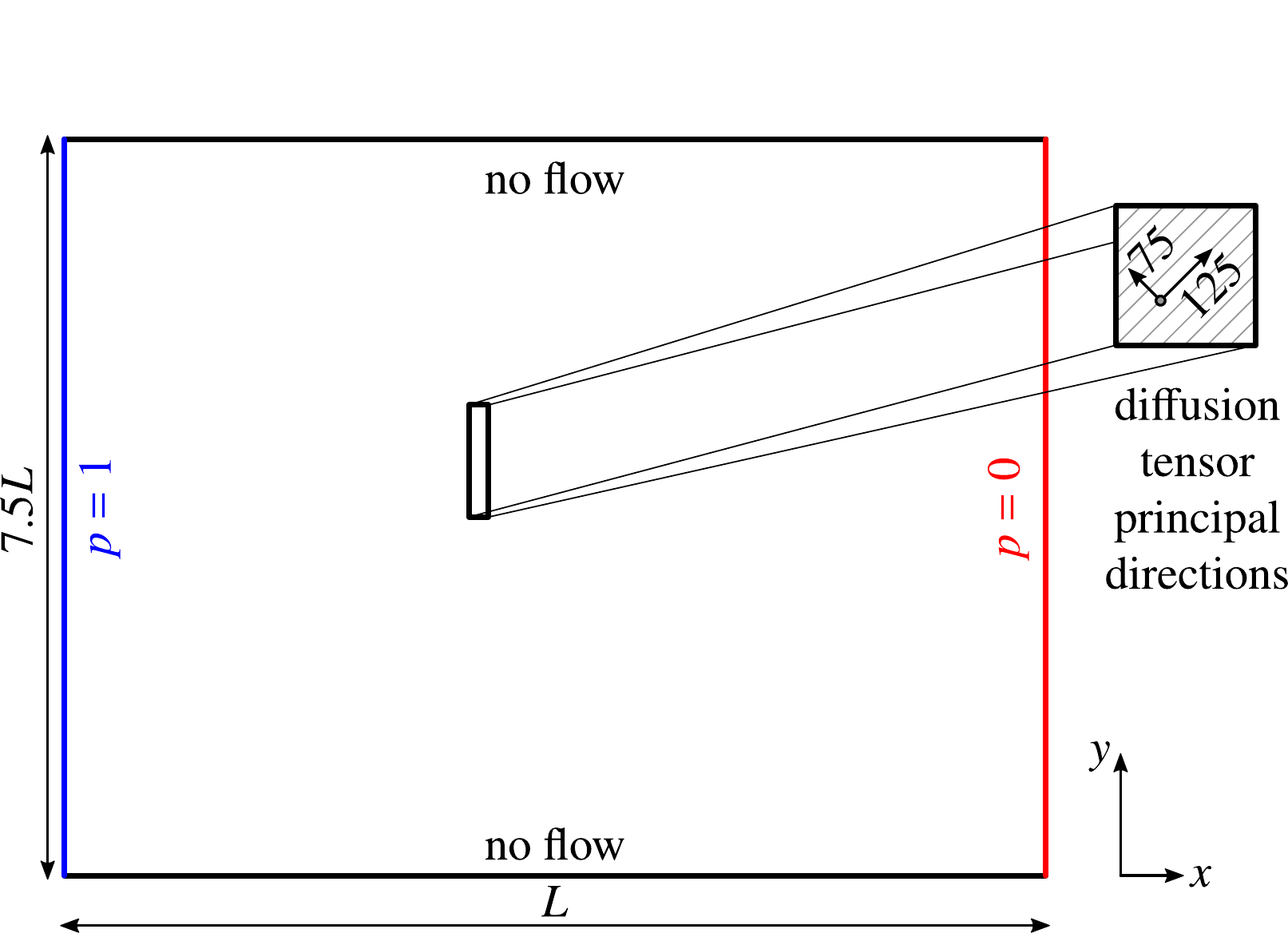}
 \caption{}
 \end{subfigure}
 \hfill
 \begin{subfigure}{0.45\textwidth}
 \includegraphics[width=\linewidth]{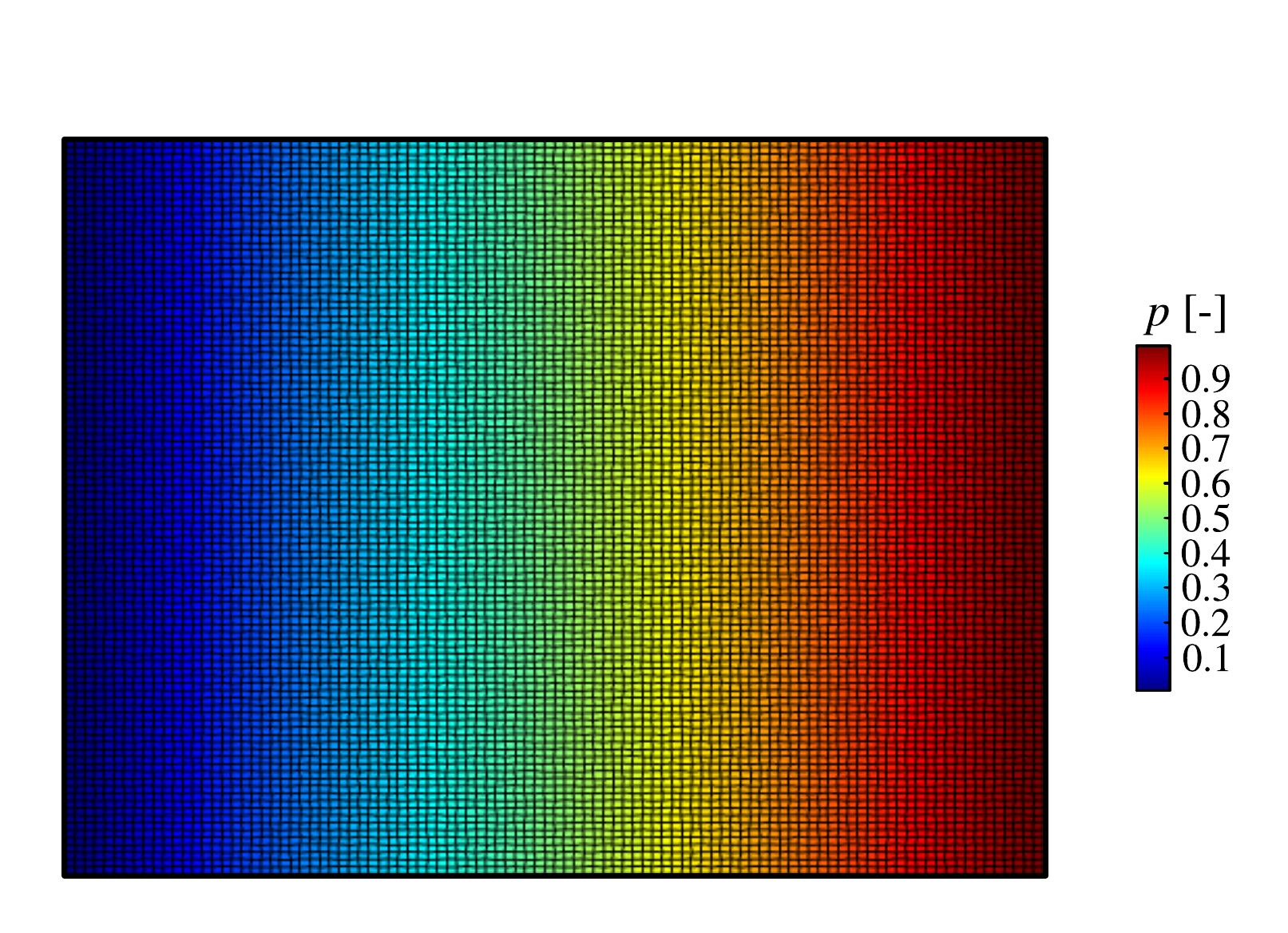}
 \caption{}
 \end{subfigure}
 \hfill\null
\caption{Grid (a) and reference solution (b) for the 2D example with a MPFA-discretized pressure equation. Note that the dimensions of the grid have been scaled for plotting visibility.}
\label{fig:mpfa_problem}
\end{figure}

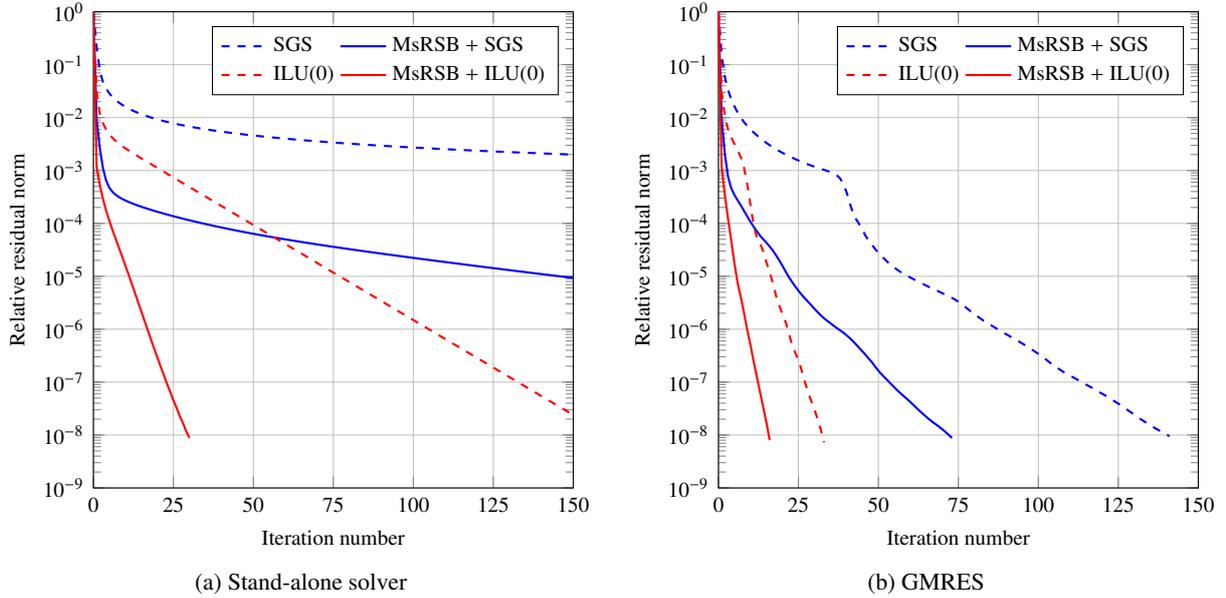
\begin{figure} [htbp]
\centering
\begin{subfigure}{0.48\textwidth}
  \begin{tikzpicture}
    \begin{semilogyaxis}[ width  = \linewidth,
                          height = 1.\linewidth,
                          grid = major,
                          major grid style={very thin,draw=gray!50},
                          xmin=0,xmax=150,
                          ymin=1e-9,ymax=1,
                          xlabel={Iteration number},
                          ylabel={Relative residual norm},
                          xtick={0,25,...,150},
                          ytick={1e-9,1e-8,1e-7,1e-6,1e-5,1e-4,1e-3,1e-2,1e-1,1},
                          ylabel near ticks,
                          xlabel near ticks,
                          tick label style={font=\footnotesize},
                          label style={font=\footnotesize},
                          legend style={font=\footnotesize},
                          legend cell align={left},
                          legend columns=2,
                          legend pos=north east]


      \addplot [dashed, blue, line width=0.8pt]
               table[x index=0, y index=1, col sep=comma]
               {2d_mpfa_RICHARDSON_smoother_SGSx1.csv};
      \addlegendentry{SGS};      
      \addplot [solid, blue, line width=0.8pt]
               table[x index=0, y index=1, col sep=comma]
               {2d_mpfa_RICHARDSON_multiscale_SGSx1.csv}; 
      \addlegendentry{MsRSB + SGS};
      
      \addplot [dashed, red, line width=0.8pt]
               table[x index=0, y index=1, col sep=comma]
               {2d_mpfa_RICHARDSON_smoother_ILU0x1.csv};
      \addlegendentry{ILU(0)};  
      \addplot [solid, red, line width=0.8pt]
               table[x index=0, y index=1, col sep=comma]
               {2d_mpfa_RICHARDSON_multiscale_ILU0x1.csv}; 
      \addlegendentry{MsRSB + ILU(0)};

    \end{semilogyaxis}
  \end{tikzpicture}
  \caption{Stand-alone solver}
\end{subfigure}
\hfill
\begin{subfigure}{0.48\textwidth}
  \begin{tikzpicture}
    \begin{semilogyaxis}[ width  = \linewidth,
                          height = 1.\linewidth,
                          grid = major,
                          major grid style={very thin,draw=gray!50},
                          xmin=0,xmax=150,
                          ymin=1e-9,ymax=1,
                          xlabel={Iteration number},
                          ylabel={Relative residual norm},
                          xtick={0,25,...,150},
                          ytick={1e-9,1e-8,1e-7,1e-6,1e-5,1e-4,1e-3,1e-2,1e-1,1},
                          ylabel near ticks,
                          xlabel near ticks,
                          tick label style={font=\footnotesize},
                          label style={font=\footnotesize},
                          legend style={font=\footnotesize},
                          legend cell align={left},
                          legend columns=2,
                          legend pos=north east]


      \addplot [dashed, blue, line width=0.8pt]
               table[x index=0, y index=1, col sep=comma]
               {2d_mpfa_GMRES_right_smoother_SGSx1.csv};
      \addlegendentry{SGS};      
      \addplot [solid, blue, line width=0.8pt]
               table[x index=0, y index=1, col sep=comma]
               {2d_mpfa_GMRES_right_multiscale_SGSx1.csv}; 
      \addlegendentry{MsRSB + SGS};
      
      \addplot [dashed, red, line width=0.8pt]
               table[x index=0, y index=1, col sep=comma]
               {2d_mpfa_GMRES_right_smoother_ILU0x1.csv};
      \addlegendentry{ILU(0)};  
      \addplot [solid, red, line width=0.8pt]
               table[x index=0, y index=1, col sep=comma]
               {2d_mpfa_GMRES_right_multiscale_ILU0x1.csv}; 
      \addlegendentry{MsRSB + ILU(0)};

    \end{semilogyaxis}
  \end{tikzpicture}
  \caption{GMRES}
\end{subfigure}
\hfill\null

\caption{ Stand-alone iterative (a) and GMRES-accelerated iterative performance (b) for the 2D MPFA example. }
\label{fig:mpfa_results}
\end{figure}

The example problem is shown in Figure~\ref{fig:mpfa_problem} and consists of a 100 by 100 structured grid discretizing a rectangular domain of 20 by 150 meters. We impose a unit pressure drop from $x=0$m to $x=20$m. Interior vertices are perturbed by a factor $0.2 \Psi \Delta x$ where $\Psi \in [-\frac{1}{2}, \frac{1}{2}]$ is a uniformly random variable with expected value 0.
The tensor $\tensorTwo{\Lambda}$ has a diagonal value of $\lambda_{xx} = \lambda_{yy} = 100$ md$\cdot$cP\textsuperscript{-1} with off-diagonal elements  $\lambda_{xy} = \lambda_{yx} = 25$ md$\cdot$cP\textsuperscript{-1}. The combined effect of the non-orthogonal grid and the strength of the off-diagonal permeability means that directly applying the MsRSB iterative process to the system matrix results in rapid divergence of the basis functions. 
With the proposed regularization, however, MsRSB can be employed and satisfactory convergence rates are obtained for both Richardson and GMRES-accelerated iterations. We partition the domain into 400 coarse blocks, each comprised of a 5 by 5 segment of fine cells and solve the MPFA system to a tolerance of $10^{-8}$. The results are displayed in Figure~\ref{fig:mpfa_results}, where Symmetric Gauss-Seidel (SGS) or ILU(0) is used as the second stage of the preconditioner (post-smoothing). No pre-smoother is applied here. We observe a clear improvement to convergence rates indicating that the basis functions successfully capture the local features of the system and resolve low-frequency errors. As a non-accelerated stand-alone solver, i.e. Richardson, MsRSB+ILU(0) converges in 30 iterations while ILU(0) fails to converge in 150 iterations. The set-up with the cheaper but less effective smoother SGS fails to converge in both cases although the multiscale stage again leads to a higher convergence rate. The performance difference between the single- and two-level solvers is reduced when GMRES is used to accelerate the solution process, nonetheless the multiscale solvers still only require half as many iterations as smoothers alone. 

\subsection{MsRSB for an MPFA-O discretization: 3D Field Test Case}
\begin{figure} [htbp]
\begin{subfigure}{0.45\textwidth}
\includegraphics[width=\linewidth]{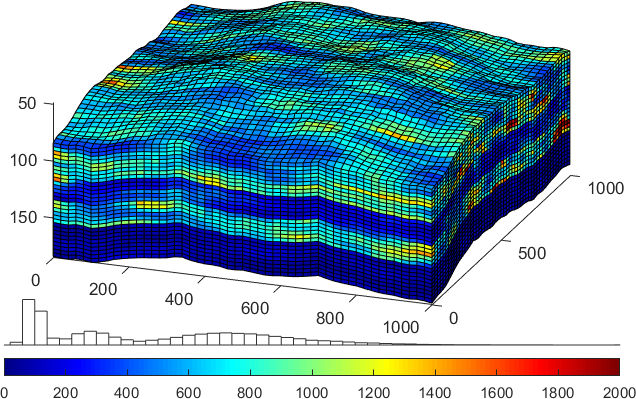}
\caption{Permeability and grid}
\label{fig:mpfa_field_perm}
\end{subfigure}
\begin{subfigure}{0.45\textwidth}
\includegraphics[width=\linewidth]{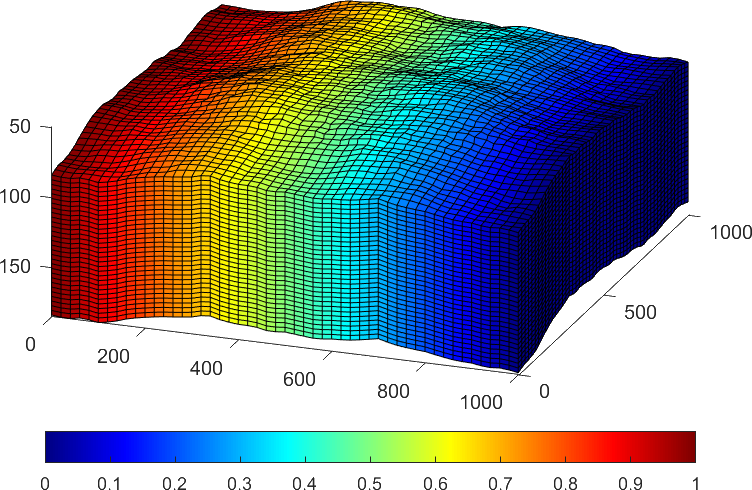}
\caption{Reference solution}
\label{fig:mpfa_field_pressure}
\end{subfigure} \\
\begin{subfigure}{0.45\textwidth}
\includegraphics[width=\linewidth]{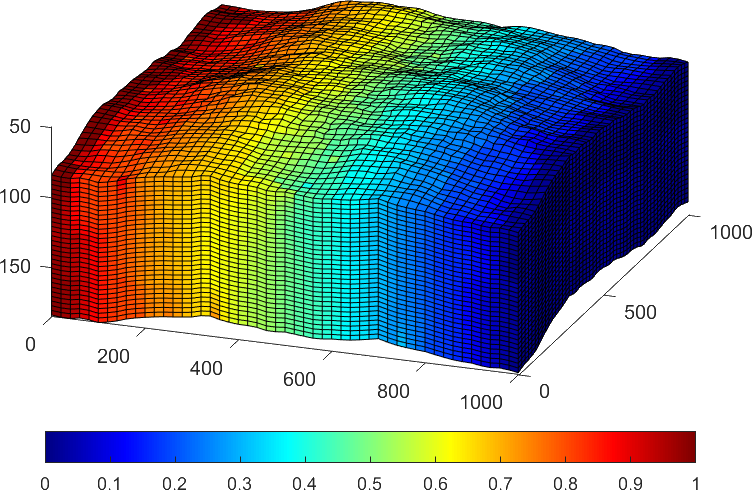}
\caption{Initial multiscale solution}
\label{fig:mpfa_field_pressure_ms}
\end{subfigure}
\begin{subfigure}{0.45\textwidth}
\includegraphics[width=\linewidth]{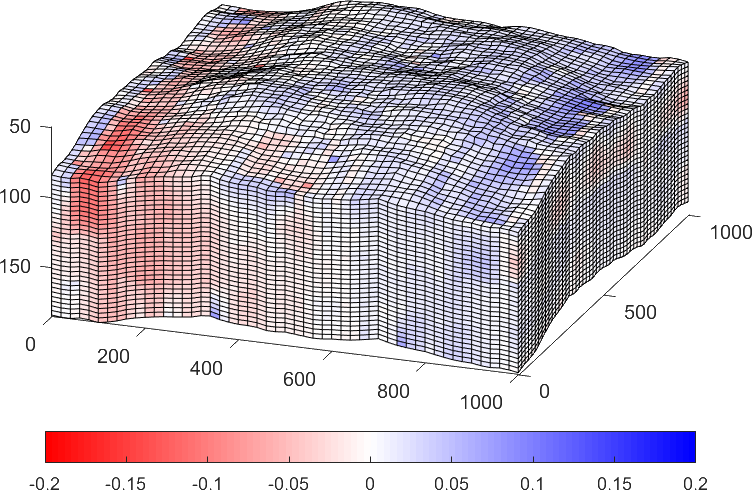}
\caption{Error in initial multiscale}
\label{fig:mpfa_field_error_ms}
\end{subfigure}
\caption{The permeability, grid (a) and the reference solution (b) for the 3D MPFA example together with the initial, uniterated multiscale solution (c) and the corresponding error relative to the fine-scale (d).}
\label{fig:mpfa_field_setup}
\end{figure}

Next a somewhat more realistic conceptual problem is considered. The test case consists of a 50 by 50 structured grid, with 30 layers in the vertical direction. The physical domain has a horizontal extent of 1000 by 1000 meters, with a vertical thickness of 100 meters. Similar to the 2D case, the vertices of the grid are perturbed to create a rough grid. Additionally, the top surface has a varying topography. The model has five different regions of a log-normally distributed diffusion tensor, with mean values of 700, 1000, 300, 800 and 100 md$\cdot$cP\textsuperscript{-1}, respectively as shown in Figure~\ref{fig:mpfa_field_perm}.

We consider a logically structured coarse mesh with block sizes of 5 by 5 by 5 fine cells, resulting in a total of 600 coarse blocks to partition the fine-grid with 75,000 cells. We apply a simple boundary condition resulting in flow from $x = 0$ to $x=1000$ meters. In engineering applications, flow would typically be driven by wells, but our goal here is to produce flow over the entire domain to verify our implementation. The reference solution is plotted in Figure~\ref{fig:mpfa_field_pressure}.

\begin{figure} [htbp]
\centering
\begin{subfigure}{0.48\textwidth}
  \begin{tikzpicture}
    \begin{semilogyaxis}[ width  = \linewidth,
                          height = 1.\linewidth,
                          grid = major,
                          major grid style={very thin,draw=gray!50},
                          xmin=0,xmax=150,
                          ymin=1e-9,ymax=1,
                          xlabel={Iteration number},
                          ylabel={Relative residual norm},
                          xtick={0,25,...,150},
                          ytick={1e-9,1e-8,1e-7,1e-6,1e-5,1e-4,1e-3,1e-2,1e-1,1},
                          ylabel near ticks,
                          xlabel near ticks,
                          tick label style={font=\footnotesize},
                          label style={font=\footnotesize},
                          legend style={font=\footnotesize},
                          legend cell align={left},
                          legend columns=2,
                          legend pos=south east]

      \addplot [dashed, blue, line width=0.8pt]
               table[x index=0, y index=1, col sep=comma]
               {3d_mpfa_RICHARDSON_smoother_SGSx1.csv};
      \addlegendentry{SGS};      
      \addplot [solid, blue, line width=0.8pt]
               table[x index=0, y index=1, col sep=comma]
               {3d_mpfa_RICHARDSON_multiscale_SGSx1.csv}; 
      \addlegendentry{MsRSB + SGS};
      
      \addplot [dashed, red, line width=0.8pt]
               table[x index=0, y index=1, col sep=comma]
               {3d_mpfa_RICHARDSON_smoother_ILU0x1.csv};
      \addlegendentry{ILU(0)};  
      \addplot [solid, red, line width=0.8pt]
               table[x index=0, y index=1, col sep=comma]
               {3d_mpfa_RICHARDSON_multiscale_ILU0x1.csv}; 
      \addlegendentry{MsRSB + ILU(0)};

    \end{semilogyaxis}
  \end{tikzpicture}
  \caption{Stand-alone solver}
\end{subfigure}
\hfill
\begin{subfigure}{0.48\textwidth}
  \begin{tikzpicture}
    \begin{semilogyaxis}[ width  = \linewidth,
                          height = 1.\linewidth,
                          grid = major,
                          major grid style={very thin,draw=gray!50},
                          xmin=0,xmax=150,
                          ymin=1e-9,ymax=1,
                          xlabel={Iteration number},
                          ylabel={Relative residual norm},
                          xtick={0,25,...,150},
                          ytick={1e-9,1e-8,1e-7,1e-6,1e-5,1e-4,1e-3,1e-2,1e-1,1},
                          ylabel near ticks,
                          xlabel near ticks,
                          tick label style={font=\footnotesize},
                          label style={font=\footnotesize},
                          legend style={font=\footnotesize},
                          legend cell align={left},
                          legend columns=2,
                          legend pos=north east]

      \addplot [dashed, blue, line width=0.8pt]
               table[x index=0, y index=1, col sep=comma]
               {3d_mpfa_GMRES_right_smoother_SGSx1.csv};
      \addlegendentry{SGS};      
      \addplot [solid, blue, line width=0.8pt]
               table[x index=0, y index=1, col sep=comma]
               {3d_mpfa_GMRES_right_multiscale_SGSx1.csv}; 
      \addlegendentry{MsRSB + SGS};
      
      \addplot [dashed, red, line width=0.8pt]
               table[x index=0, y index=1, col sep=comma]
               {3d_mpfa_GMRES_right_smoother_ILU0x1.csv};
      \addlegendentry{ILU(0)};  
      \addplot [solid, red, line width=0.8pt]
               table[x index=0, y index=1, col sep=comma]
               {3d_mpfa_GMRES_right_multiscale_ILU0x1.csv}; 
      \addlegendentry{MsRSB + ILU(0)};

    \end{semilogyaxis}
  \end{tikzpicture}
  \caption{GMRES}
\end{subfigure}
\hfill\null

\caption{ Stand-alone iterative (a) and GMRES-accelerated iterative performance (b) for the 3D MPFA example. }
\label{fig:mpfa_field_results}
\end{figure}
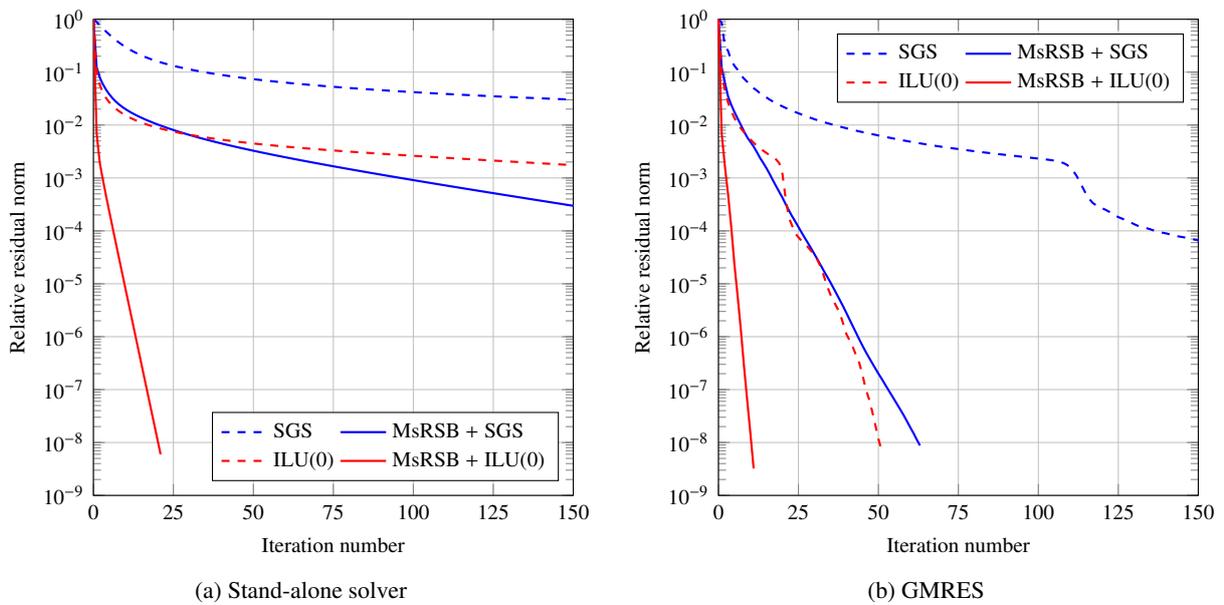

The convergence rates of MsRSB+ILU(0), ILU(0), MsRSB+SGS and SGS are compared in Figure~\ref{fig:mpfa_field_results} with and without Krylov acceleration. We observe that the benefits of the multiscale stage are more significant than in the 2D example. This is likely due to the denser stencil in 3D which leads to a fundamentally more difficult linear system, even for an equivalent number of degrees of freedom. We point out that the multiscale solver with SGS smoothing exhibits robust convergence with GMRES acceleration while the smoother-only setup stagnates. We also observe that the performance of the multiscale solver with ILU(0) is excellent, using 21 and 11 iterations without and with GMRES, respectively. We can also examine the initial multiscale solution, that is, the solution without any application of ILU(0) as an approximate solver. The solution is shown in Figure~\ref{fig:mpfa_field_pressure_ms} where we observe that the general solution is accurately captured, with minor artifacts near the boundary of the domain. The solution in the interior of the domain is largely unaffected, which is reflected in the sum of the error: $\| \mathbf{p} - \mathbf{p}_{ms} \|_1 / \| \mathbf{p} \|_1 = 0.0385$ with a maximum cell-wise error of 0.1883. These values are in line with similar MsRSB-TPFA examples for 3D models where flow is driven by boundary conditions.

\section{MsRSB for non-M matrices: FE simulation of linear elastic geomechanics} \label{sec:geomechanics}
In this section, we investigate the extension of the enhanced MsRSB preconditioner to linear systems arising from geomechanical problems assuming linear elastic behavior. By means of this application, we extend the method's applicability to problems with vector unknowns.
The linear elastostatics governing equations and their FE discretization using a classical displacement formulation \cite{Hug00} are reviewed in \ref{app:model_elastostatics} and \ref{app:model_elasticity_FE}, respectively.
To compute the nodal displacement discrete solution, the FE method requires the solution of linear systems characterized by a fine scale stiffness matrix $A$ that is SPD but not an M-matrix.
Therefore, in the construction of the prolongation operator, a na\"{i}ve implementation of MsRSB would exhibit the same convergence issues as discussed for the pressure equation in Section \ref{sec:fix_non_M-matrix}.

To apply the enhanced MsRSB preconditioner, we observe that if the displacement degrees of freedom are ordered based on each coordinate direction, the stiffness matrix $A$ possesses a block structure, namely:

\begin{align}
	A &=
    \begin{bmatrix}
  	  A_{xx} & A_{xy} & A_{xz} \\
  	  A_{yx} & A_{yy} & A_{yz} \\
  	  A_{zx} & A_{zy} & A_{zz} 
	\end{bmatrix},
	\label{eq:blk_stiff}
\end{align} 

\noindent
which reflects the full coupling between $x$, $y$, and $z$ components of displacements.
For preconditioning purposes, the complete matrix is often replaced with a sparser block diagonal approximation

\begin{align}
	A^{\textsc{(sdc)}} &=
    \begin{bmatrix}
  	  A_{xx} &        &        \\
  	         & A_{yy} &        \\
  	         &        & A_{zz} 
	\end{bmatrix},
  \label{eq:blk_stiff_SDC}
\end{align} 

\noindent
namely the \textit{separate displacement component} (SDC) approximation proposed in \cite{AxeGus78}.
The motivation behind the SDC approximation is that, using Korn's inequality, one can show that $\Mat{A}^\text{(SDC)}$ is spectrally equivalent to $\Mat{A}$ \cite{Bla94,GusLin98}.
Note that this approximation breaks down in the incompressible elasticity limit, i.e. Poisson ratio $\nu \to 0.5$.

As each diagonal block in \eqref{eq:blk_stiff} corresponds to the finite element discretization of an anisotropic diffusion operator for the corresponding displacement component---see Remark \ref{rem:SDC}---, enhanced MsRSB can be readily applied to \eqref{eq:blk_stiff_SDC} to obtain an approximate solution or preconditioner for the fine-scale matrix.
In two- (or three-) dimensional FE-based elasticity simulation, two (or three, respectively) basis functions are associated with each coarse node.
Computing such basis functions using $\Mat{A}^{\textsc{(sdc)}}$ implies that the fine-scale displacement field in each coordinate direction is expressed as a linear combination of coarse nodal displacement in the corresponding direction only. Hence, the prolongation operator will be block diagonal.
Nevertheless, we emphasize that the coarse operator, $\Mat{A}_c$, computed with eq. \eqref{eq:lin_system_coarse}, will still capture full coupling influences between $x$-, $y$-, and $z$-displacement components when using enhanced MsRSB basis functions.
Note that the block diagonal property of the prolongation operator represents a major difference with a classic MSFE approach \cite{MultiscaleFEM_Castelleto2017}. At the cost of a much denser prolongation operator, those approaches account for additional coupling between all displacement solution directions and basis functions directions.

\begin{figure} [htbp]
\begin{subfigure}[t]{0.22\textwidth}
  \centerline{\includegraphics[width=\linewidth]{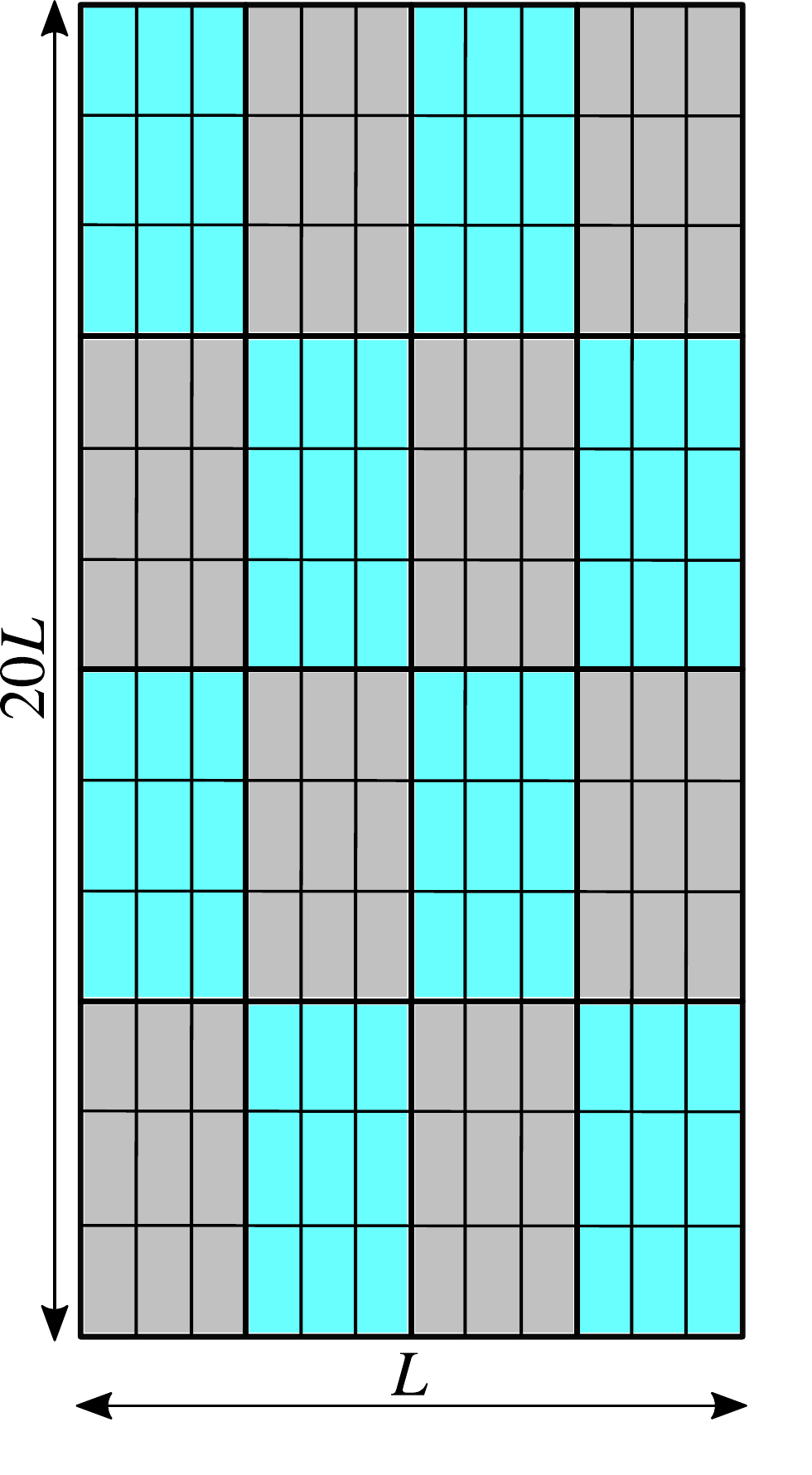}}
  \caption{\label{fig:FE_demo_a}}
\end{subfigure}
\hfill
\begin{subfigure}[t]{0.22\textwidth}
  \centerline{\includegraphics[width=\linewidth]{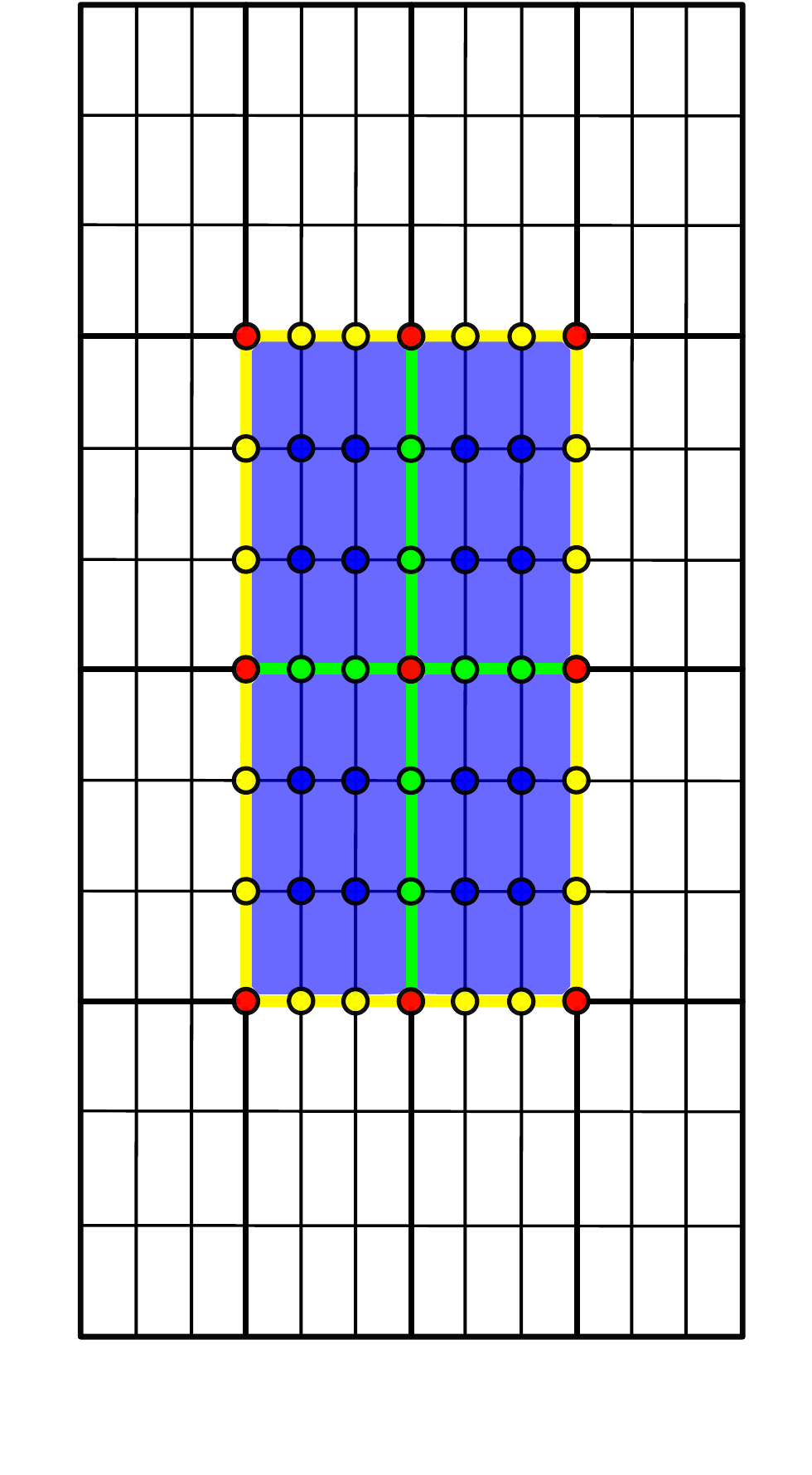}}
  \caption{\label{fig:FE_demo_b}}
\end{subfigure}
\hfill
\begin{subfigure}[t]{0.22\textwidth}
  \centerline{\includegraphics[width=\linewidth]{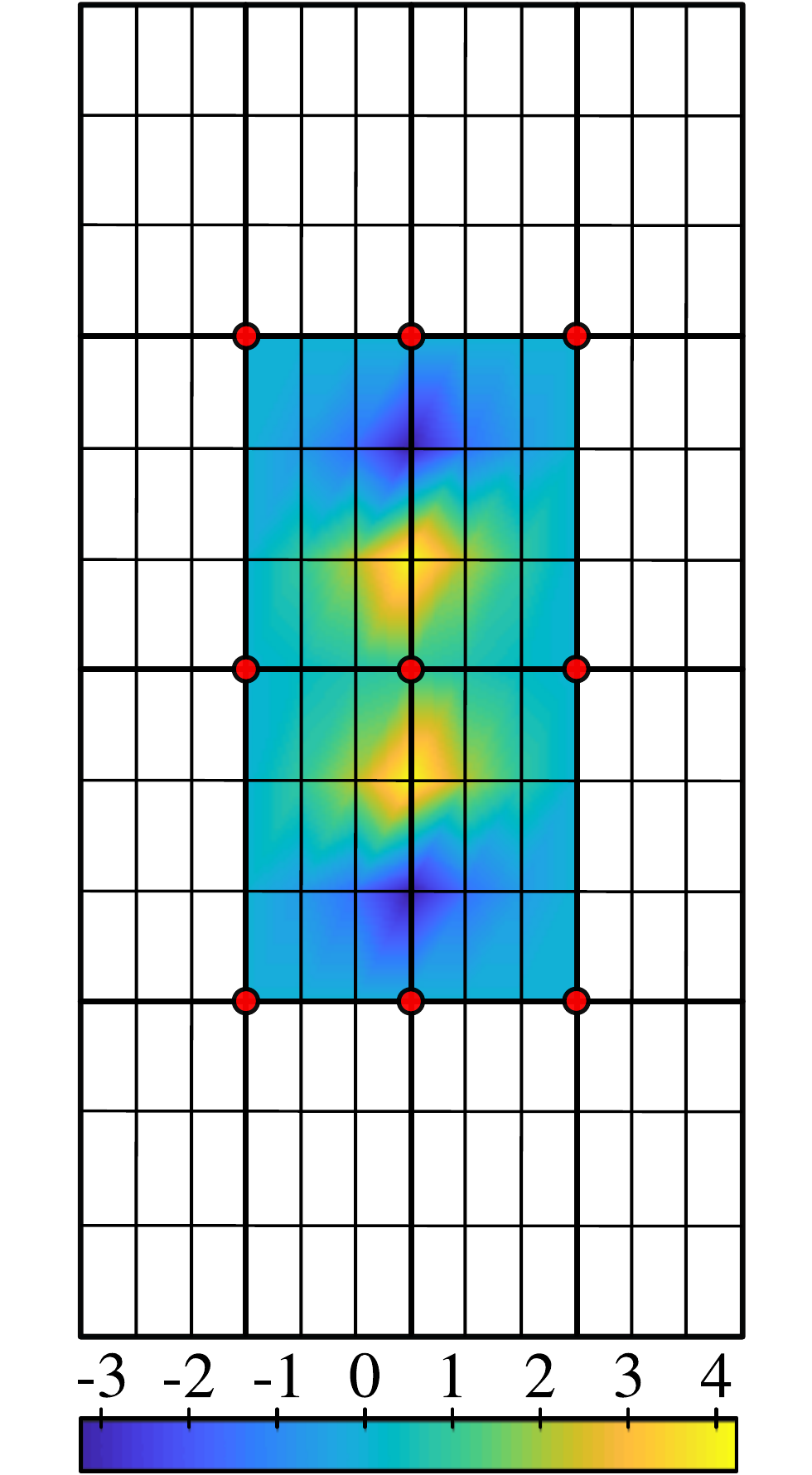}}
  \caption{\label{fig:FE_demo_c}}
\end{subfigure}
\hfill
\begin{subfigure}[t]{0.22\textwidth}
  \centerline{\includegraphics[width=\linewidth]{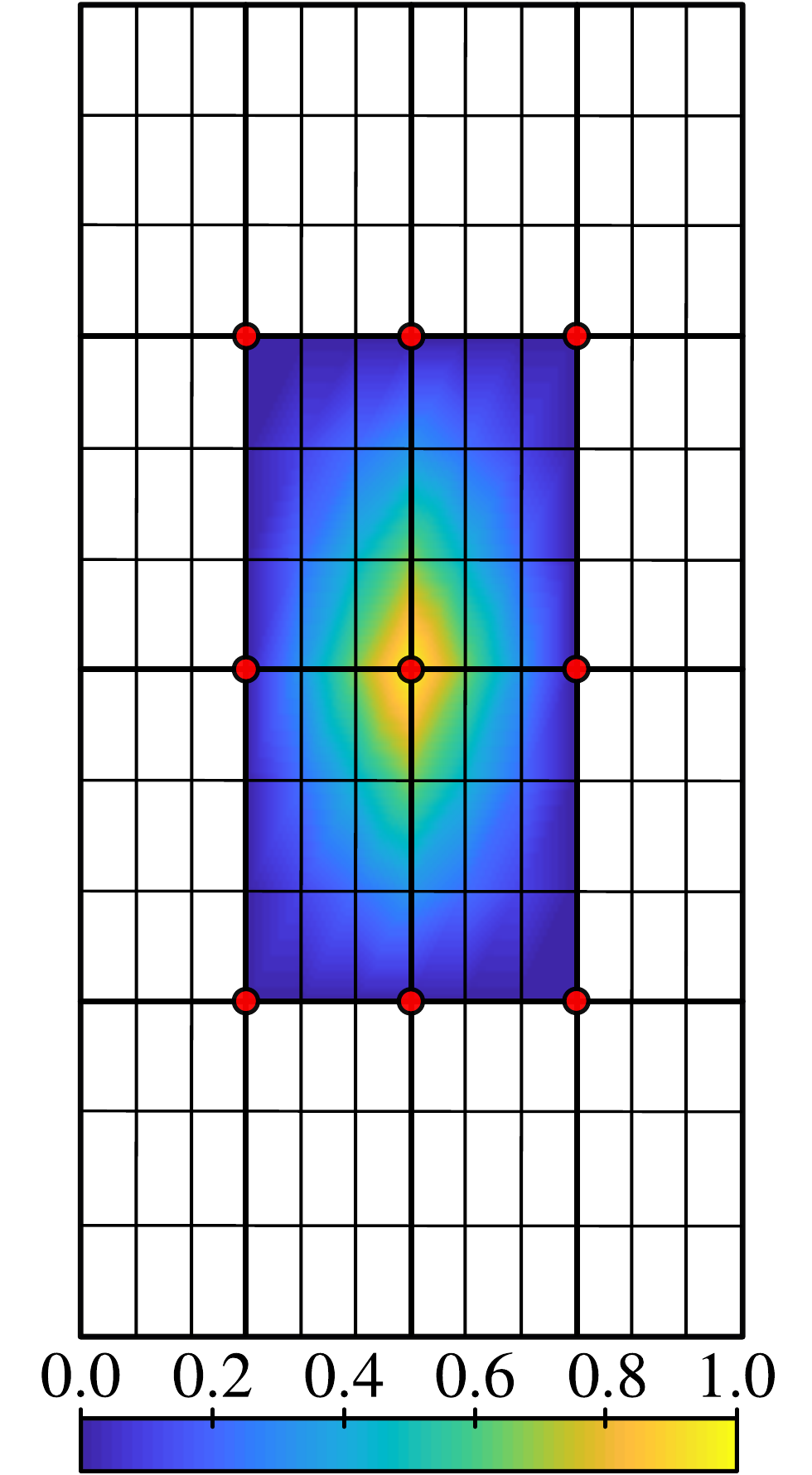}}
  \caption{\label{fig:FE_demo_d}}
\end{subfigure}
\caption{\label{fig:bfSimpleTestGeomechAdapted} Basis function relative to the displacement degree of freedom in the $x$-direction associated to the innermost coarse node of a regular 4$\times$4 partition: (a) primal grid; (b) coarse (\protect\tikz \protect\draw[black,fill=red] (0,0) circle (.5ex);), support boundary (\protect\tikz \protect\draw[black,fill=yellow] (0,0) circle (.5ex);), edge (\protect\tikz \protect\draw[black,fill=green] (0,0) circle (.5ex);) and internal (\protect\tikz \protect\draw[black,fill=blue] (0,0) circle (.5ex);) nodes; (c) MsRSB basis function using the original fine-scale system; and (d) enhanced MsRSB basis function using the filtered fine-scale system. The convergence criterion is $e_{it} < 10^{-3}$, where $e_{it}$ is the maximum value in the smoother update on the internal nodes (see algorithm \ref{alg:EnhancedMsRSB}). However as the case with the original linear system diverges, the basis function obtained after 20 iterations is plotted in (c).}
\end{figure}

This extension to the enhanced MsRSB method is tested on a simple homogeneous problem defined in terms of dimensionless quantities with unit Lam\'e parameters, i.e. setting $E$=1 and $\nu$=0.25.
The fine grid is chosen to be $12\times12$ Cartesian.
The coarse grid has 5 coarse nodes in each direction.
Finally, the grid is rescaled by a factor of $20$ in the $y$-dimension. The high aspect ratio is chosen to induce strong non-M matrix properties of the resulting linear system. 
Figure \ref{fig:FE_demo_a} depicts the coarse and fine grid of the test case on the $x$-$y$ plane.
To assess the robustness of the enhanced MsRSB method including the extension to vector physics, a basis function obtained using a na\"ive implementation of MsRSB is compared to the same basis function obtained with the enhanced method. 
Here, a na\"ive implementation of MsRSB is a straightforward application of the original method for M matrices to the linear system. 
Figure \ref{fig:bfSimpleTestGeomechAdapted} depicts the obtained basis functions for a given coarse node. It is obvious from the plots that the original MsRSB produces to diverging basis functions whereas the enhanced method recovers the expected bi-linear interpolators of a Cartesian homogeneous problem.
The results reemphasize the findings of section \ref{sec:fix_non_M-matrix}.

\subsection{MsRSB for Geomechanics: A 2D Heterogeneous Test Case}

\begin{figure}[htbp]
    \centering
    \includegraphics[width=\textwidth]{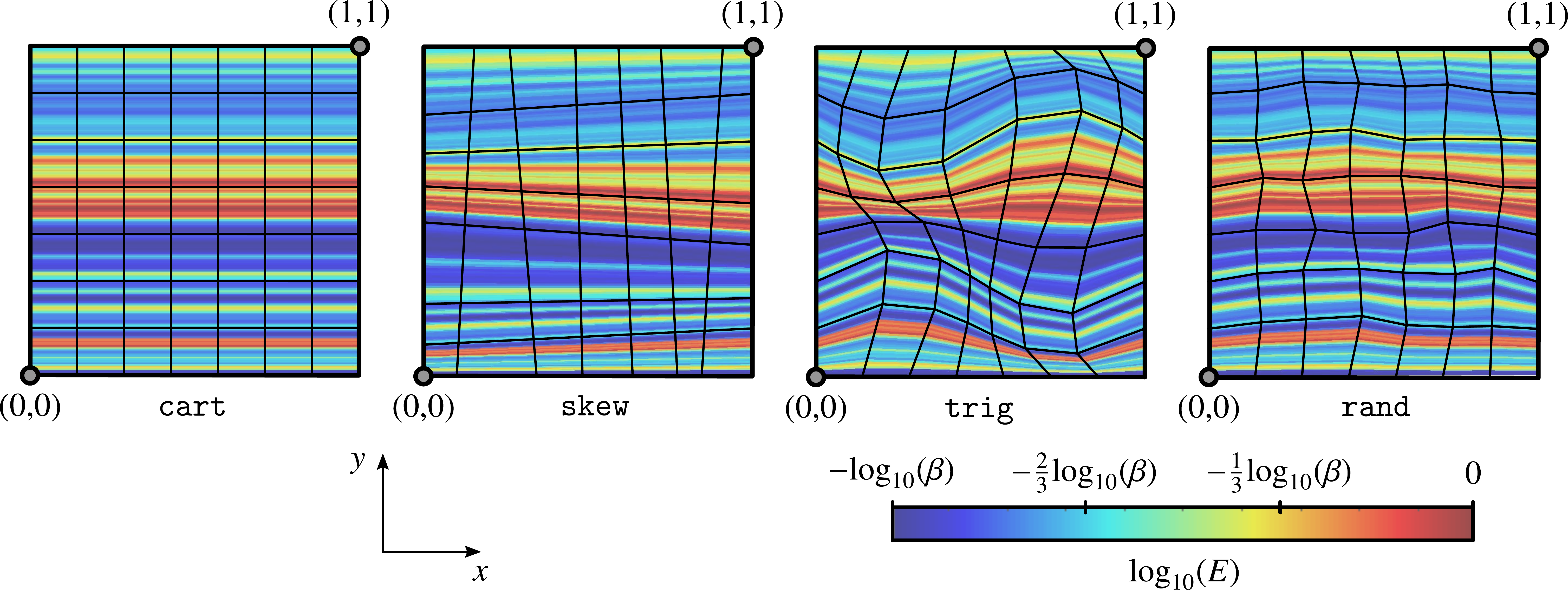}
    \caption{Structured 2D heterogeneous case \cite{MultiscaleFEM_Castelleto2017}: Young's modulus distribution and computational mesh setup.}
    \label{fig:2D_structured_young}
\end{figure}

To compare the enhanced MsRSB method to an existing multiscale preconditioner for geomechanics \cite{MultiscaleFEM_Castelleto2017}, we investigate a 2D test case consisting of an elastic isotropic domain with a heterogeneous (layered) distribution of Young's modulus, shown in Figure~\ref{fig:2D_structured_young}. Four mesh families (\texttt{cart}, \texttt{skew}, \texttt{trig} and \texttt{rand}) are considered, each with a different type of geometric distortion applied to the coarse elements. In each case a $224\times224$ fine-scale mesh and $7\times7$ coarse mesh are defined on the domain. The reader is referred to section 4.1 of \cite{MultiscaleFEM_Castelleto2017} for a more detailed description of the test case setup.   Two Krylov solver setups are tested.   The first variant is Preconditioned Conjugate Gradient (PCG), with a symmetric preconditioning operator constructed by pre-- and post--smoothing the multiscale (MsRSB) operator with no-fill Incomplete Cholesky factorization (IC(0)).   The second version employs Biconjugate Gradient Stabilized (BiCGStab) with a two--stage preconditioning scheme, which consists of the multiscale operator as the global step followed by no--fill incomplete LU factorization (ILU(0)) as post-smoothing with no pre-smoother.   This setup is identical to the one used in \cite{MultiscaleFEM_Castelleto2017}, except that MSFE is replaced by enhanced MsRSB to construct the multiscale operator.  

Table~\ref{tab:2D_structured_CG_iter} summarizes the iteration counts using PCG for different mesh families and levels of contrast in material properties.  With exception of the case with extremely large coarse elements ($7\times7$, which results in each coarse element consisting of $32\times32$ fine elements), the preconditioner offers robust performance over a variety of mesh families and types of boundary conditions. In a direct comparison against an identical setup using MSFE, from \cite{MultiscaleFEM_Castelleto2017},  Table~\ref{tab:2D_structured_BICGSTAB_iter} shows similar results obtained with BiCGSTab and a single ILU(0) smoothing step. Here, the proposed MsRSB method achieves similar performance compared to MSFE, with only a modest (about 20--25\% on average) increase in the number of iterations for most cases, with a notable exception of Cartesian laterally constrained case, where MSFE basis functions provide an exact interpolation. Note that the BiCGStab solver applies the preconditioner twice at every iteration, which corresponds to two applications of both multiscale operator and smoother, whereas PCG employs a single application of the preconditioner, which involves two smoothing steps (pre-- and post--smoothing) and only one multiscale step. Nevertheless, comparable performance is observed in terms of iteration count, suggesting PCG with symmetric preconditioning should lead to a more efficient option due to its lower cost per iteration. 

\begin{table*}
    \centering
    \caption{Structured 2D heterogeneous case: number of PCG iterations to converge to $10^{-8}$ for different mesh families, coarse grid sizes and Young's modulus variations.}
    \label{tab:2D_structured_CG_iter}
    \begin{small}
    \setlength\tabcolsep{3.8pt}
    \begin{tabular}{cccccccccccccccc}
        \hline\noalign{\smallskip}
        $\beta$ & \multirow{2}{*}{\begin{tabular}{c}\# coarse\\elements\end{tabular}} & \multicolumn{4}{c}{laterally constrained} & & \multicolumn{4}{c}{laterally unconstrained} & & \multicolumn{4}{c}{simple shear} \\
        \noalign{\smallskip}\cline{3-6} \cline{8-11} \cline{13-16}\noalign{\smallskip}
        & & \texttt{cart} & \texttt{skew} & \texttt{trig} & \texttt{rand} & & \texttt{cart} & \texttt{skew} & \texttt{trig} & \texttt{rand} & & \texttt{cart} & \texttt{skew} & \texttt{trig} & \texttt{rand} \\
        \hline\noalign{\smallskip}
        1 & $56\times56$ & 5 & 7 & 7 & 6 & & 5 & 7 & 7 & 6 & & 6 & 7 & 7 & 6 \\
          & $28\times28$ & 9 & 12 & 13 & 10 & & 9 & 13 & 13 & 11 & & 12 & 15 & 15 & 13 \\
          & $14\times14$ & 16 & 21 & 22 & 20 & & 18 & 24 & 23 & 22 & & 24 & 27 & 27 & 23 \\
          & $7 \times7 $ & 37 & 40 & 46 & 40 & & 41 & 45 & 48 & 45 & & 50 & 43 & 47 & 43 \\
        \hline\noalign{\smallskip}
        2 & $56\times56$ & 5 & 6 & 7 & 6 & & 7 & 7 & 8 & 7 & & 7 & 7 & 8 & 7 \\
          & $28\times28$ & 8 & 11 & 12 & 10 & & 9 & 12 & 12 & 10 & & 13 & 15 & 17 & 14 \\
          & $14\times14$ & 17 & 22 & 25 & 21 & & 19 & 24 & 27 & 24 & & 32 & 32 & 32 & 30 \\
          & $7 \times7 $ & 47 & 44 & 52 & 48 & & 58 & 53 & 55 & 55 & & 67 & 50 & 54 & 52 \\
        \hline\noalign{\smallskip}
        3 & $56\times56$ & 5 & 7 & 10 & 8 & & 9 & 9 & 11 & 9 & & 10 & 10 & 11 & 10 \\
          & $28\times28$ & 8 & 12 & 14 & 11 & & 11 & 13 & 15 & 13 & & 16 & 17 & 19 & 16 \\
          & $14\times14$ & 21 & 28 & 42 & 30 & & 23 & 31 & 44 & 33 & & 39 & 40 & 52 & 40 \\
          & $7 \times7 $ & 50 & 54 & 66 & 60 & & 72 & 64 & 85 & 71 & & 94 & 63 & 75 & 64 \\
        \hline\noalign{\smallskip}
    \end{tabular}
    \end{small}
\end{table*}

\begin{table*}
    \centering
    \caption{Structured 2D heterogeneous case: number of BICGSTAB iterations to converge to $10^{-8}$ for different mesh families, coarse grid sizes and Young's modulus variations. }
    \label{tab:2D_structured_BICGSTAB_iter}
    \begin{small}
    \setlength\tabcolsep{3.8pt}
    \begin{tabular}{cccccccccccccccc}
        \hline\noalign{\smallskip}
        $\beta$ & \multirow{2}{*}{\begin{tabular}{c}\# coarse\\elements\end{tabular}} & \multicolumn{4}{c}{laterally constrained} & & \multicolumn{4}{c}{laterally unconstrained} & & \multicolumn{4}{c}{simple shear} \\
        \noalign{\smallskip}\cline{3-6} \cline{8-11} \cline{13-16}\noalign{\smallskip}
        & & \texttt{cart} & \texttt{skew} & \texttt{trig} & \texttt{rand} & & \texttt{cart} & \texttt{skew} & \texttt{trig} & \texttt{rand} & & \texttt{cart} & \texttt{skew} & \texttt{trig} & \texttt{rand} \\
        \hline\noalign{\smallskip}
        1 & $56\times56$ & 5 & 7 & 6 & 5 & & 4 & 7 & 6 & 5 & & 6 & 7 & 7 & 5 \\
          & $28\times28$ & 8 & 10 & 11 & 8 & & 8 & 12 & 11 & 10 & & 13 & 16 & 17 & 14 \\
          & $14\times14$ & 16 & 20 & 19 & 16 & & 18 & 21 & 22 & 19 & & 26 & 29 & 32 & 26 \\
          & $7 \times7 $ & 37 & 33 & 40 & 40 & & 41 & 38 & 43 & 43 & & 40 & 47 & 44 & 42 \\
        \hline\noalign{\smallskip}
        2 & $56\times56$ & 4 & 6 & 6 & 5 & & 5 & 6 & 6 & 6 & & 5 & 6 & 6.5 & 5 \\
          & $28\times28$ & 7 & 10 & 12 & 8 & & 9 & 10 & 12 & 9 & & 12 & 19 & 20 & 17 \\
          & $14\times14$ & 17 & 26 & 26 & 18 & & 18 & 23 & 25 & 22 & & 37 & 39 & 46 & 37 \\
          & $7 \times7 $ & 46 & 50 & 54 & 49 & & 52 & 58 & 53 & 53 & & 51 & 59 & 55 & 51 \\
        \hline\noalign{\smallskip}
        3 & $56\times56$ & 4 & 6 & 8 & 6 & & 6 & 8 & 10 & 9 & & 6 & 7 & 8 & 7 \\
          & $28\times28$ & 10 & 11 & 13 & 10 & & 11 & 12 & 17 & 10 & & 17 & 19 & 19 & 18 \\
          & $14\times14$ & 22 & 30 & 39 & 30 & & 21 & 32 & 40 & 34 & & 46 & 49 & 47 & 50 \\
          & $7 \times7 $ & 49 & 61 & 73 & 63 & & 63 & 71 & 86 & 71 & & 62 & 67 & 71 & 58 \\
        \hline\noalign{\smallskip}
    \end{tabular}
    \end{small}
\end{table*}

\subsection{MsRSB for Geomechanics: A Geological 2D Cross-Section Test Case}

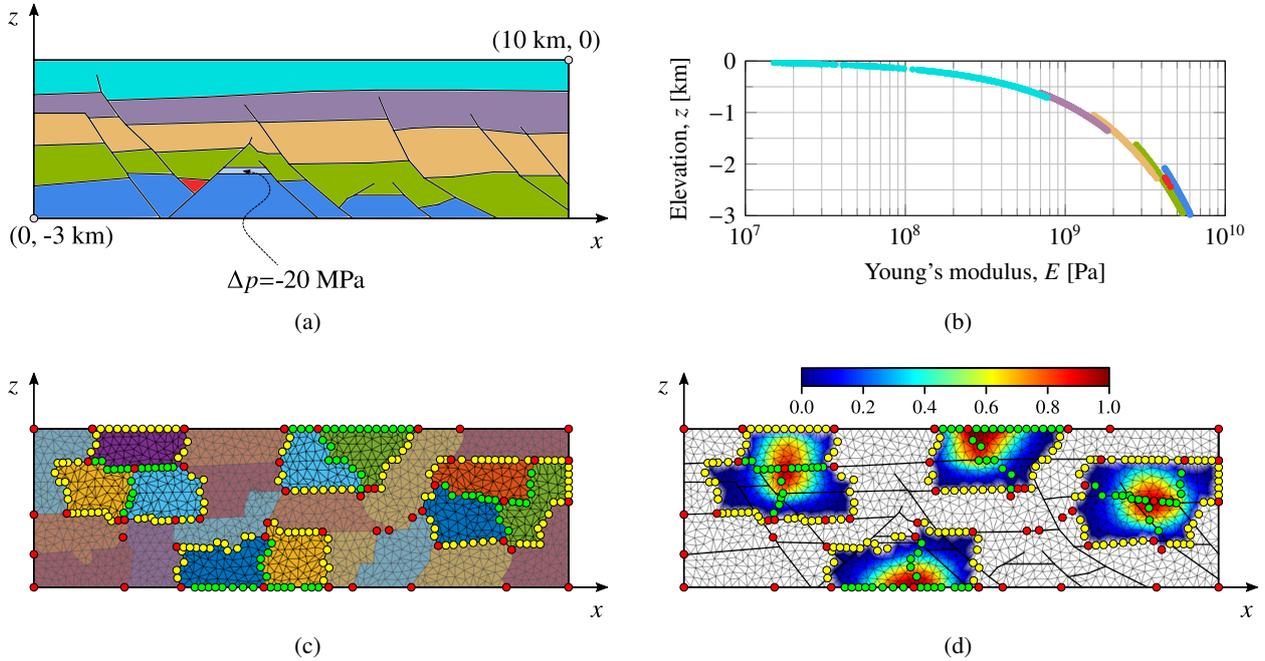
\begin{figure}[htbp]
   \centering
    \begin{subfigure}[t]{0.48\textwidth}
        \centering
        \includegraphics[align=b,width=\textwidth]{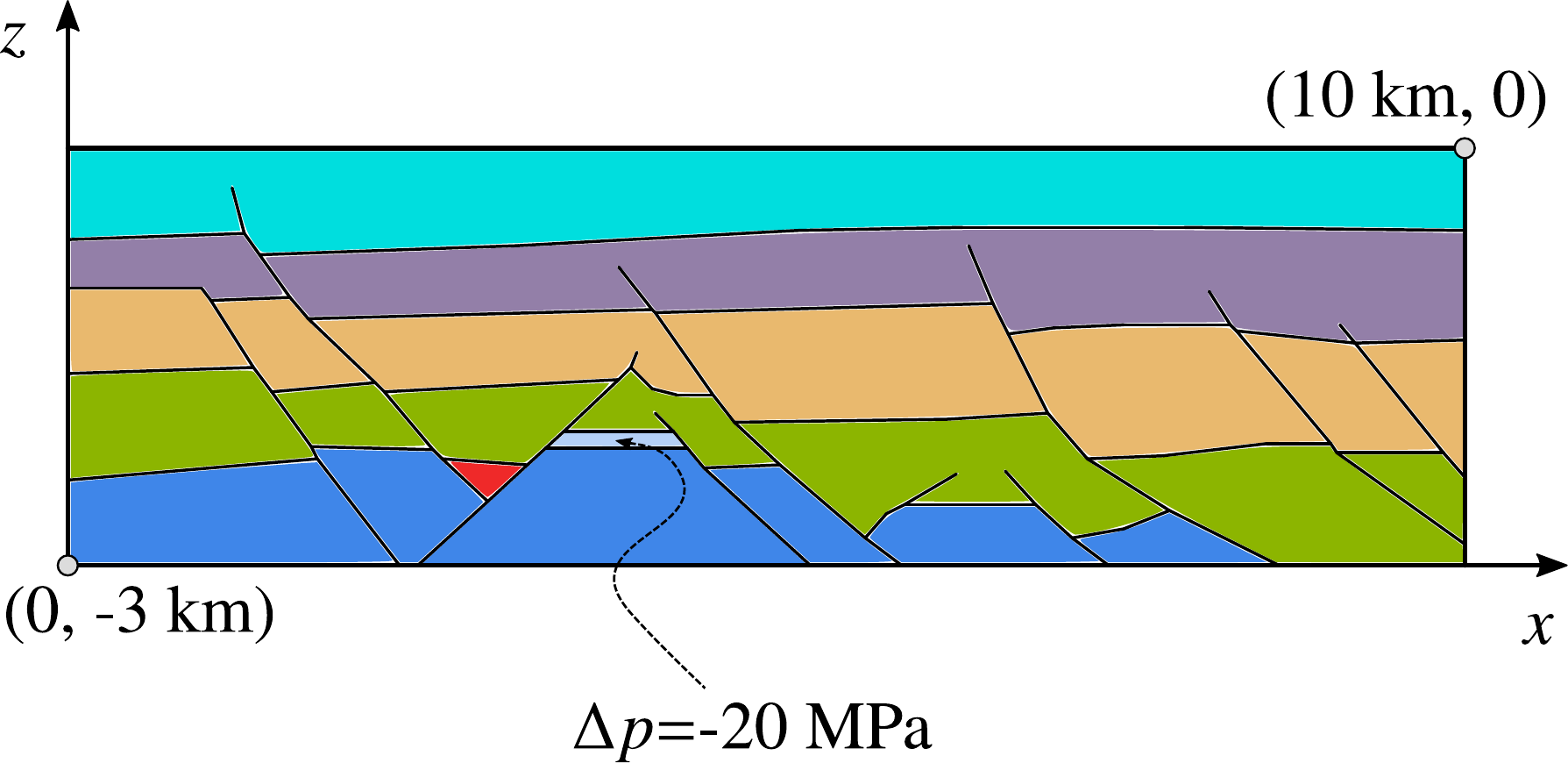}
        \caption{}
        \label{fig:cross_section_sketch}
    \end{subfigure}
    \hfill    
    \begin{subfigure}[t]{0.48\textwidth}
        \centering
 \begin{tikzpicture}
    \definecolor{layer1}{rgb}{0.247,0.524,0.912}
    \definecolor{layer2}{rgb}{0.937,0.161,0.161}
    \definecolor{layer3}{rgb}{0.550,0.710,0.000}
    \definecolor{layer4}{rgb}{0.914,0.725,0.431}
    \definecolor{layer5}{rgb}{0.678,0.498,0.659}
    \definecolor{layer6}{rgb}{0.000,0.870,0.870}
    \begin{semilogxaxis}[ width  = \linewidth,
                          height = .46\linewidth,
                          grid = both,
                          major grid style={very thin,draw=gray!50},
                          xmin=1e7,xmax=1e10,
                          ymin=-3,ymax=0,
                          xlabel={Young's modulus, $E$ [Pa]},
                          ylabel={Elevation, $z$ [km]},
                          xtick={1e7,1e8,1e9,1e10},
                          ytick={-3,-2,...,0},
                          minor y tick num = 1,
                          ylabel near ticks,
                          xlabel near ticks,
                          tick label style={font=\small},
                          label style={font=\small},
                          legend style={font=\small},
                          legend cell align={left},
                          legend columns=2,
                          legend pos=north east]

      \addplot+[layer1, only marks, mark=*, mark size =1pt,mark options={solid,fill=layer1}]
               table[x=young_Pa, y=elevation_m, col sep=comma]
               {layer_1.csv};
               
      \addplot+[layer3, only marks, mark=*, mark size =1pt,mark options={solid,fill=layer3}]
               table[x=young_Pa, y=elevation_m, col sep=comma]
               {layer_3.csv};   

      \addplot+[layer2, only marks, mark=*, mark size =1pt,mark options={solid,fill=layer2}]
               table[x=young_Pa, y=elevation_m, col sep=comma]
               {layer_2.csv};            

      \addplot+[layer4, only marks, mark=*, mark size =1pt,mark options={solid,fill=layer4}]
               table[x=young_Pa, y=elevation_m, col sep=comma]
               {layer_4.csv};

      \addplot+[layer5, only marks, mark=*, mark size =1pt,mark options={solid,fill=layer5}]
               table[x=young_Pa, y=elevation_m, col sep=comma]
               {layer_5.csv};
               
      \addplot+[layer6, only marks, mark=*, mark size =1pt,mark options={solid,fill=layer6}]
               table[x=young_Pa, y=elevation_m, col sep=comma]
               {layer_6.csv};

    \end{semilogxaxis}
  \end{tikzpicture}
        \caption{}
        \label{fig:cross_section_young_2D}
    \end{subfigure}

    \begin{subfigure}[t]{0.48\textwidth}
        \centering
        \includegraphics[align=b,width=\textwidth]{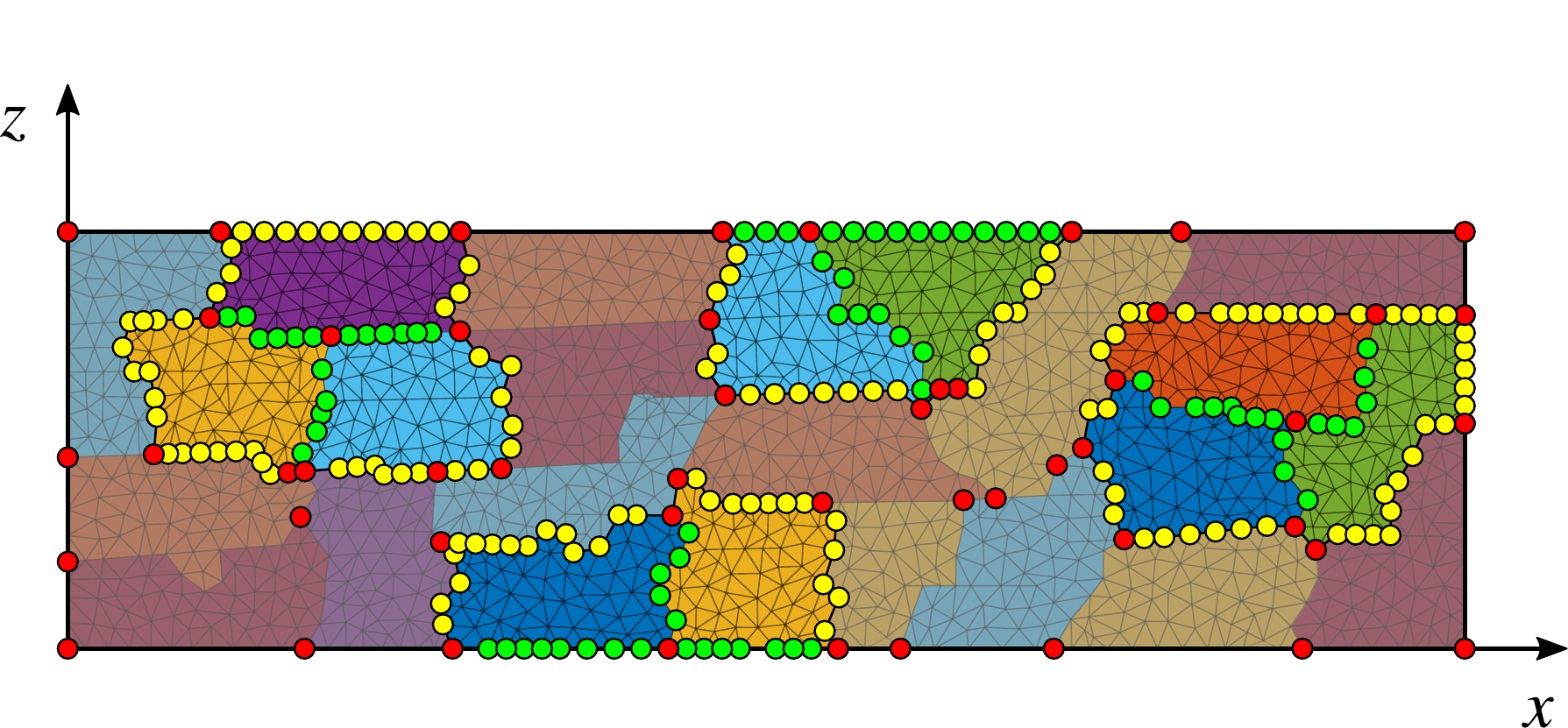}
        \caption{}
        \label{fig:cross_section_mesh}
    \end{subfigure}
    \hfill
    \begin{subfigure}[t]{0.48\textwidth}
        \centering
        \includegraphics[align=b,width=\textwidth]{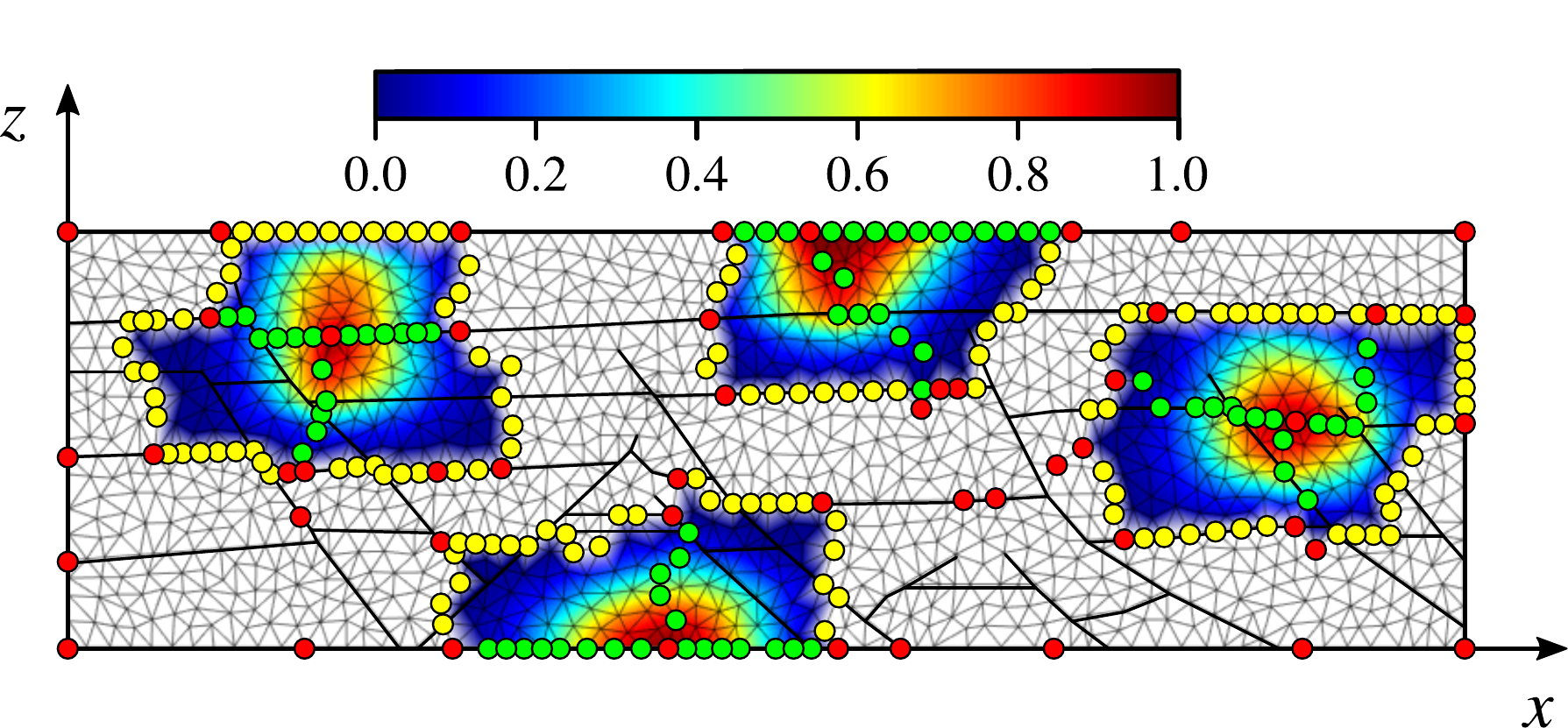}
        \caption{}
        \label{fig:crosssection_basis}
    \end{subfigure}
    \caption{Geological 2D cross--section test case: (a) physical domain, layers and faults; (b) Young's modulus distribution as a function of the elevation for each layer shown in (a); (c) fine and coarse scale grids in addition to support regions of selected basis functions (red markers denote coarse nodes; yellow markers are support boundaries; green markers are support edges); (d) examples of MsRSB displacement basis functions after 10 iterations and their support regions.}
    \label{fig:cross_section}
\end{figure}

A second geomechanics test case is designed to represent a 2D cross-section ($x$-$z$ plane) of an elastic subsurface porous medium domain, characterized by distinct geological layers and faults (shown in Figure~\ref{fig:cross_section_sketch}). A vertical distribution of Young's modulus is prescribed, based on a correlation for uniaxial compressibility developed in \cite{Bau_etal02} and recently used in \cite{MultiscaleFEM_Castelleto2017}. Specifically, the medium vertical compressibility is computed as
\begin{linenomath}
\begin{align}
    c_M = 0.01241 \,|\sigma'_z|^{-1.1342} \label{eq:compressibility2Dtestcase}
\end{align}
\end{linenomath}
where
\begin{linenomath}
\begin{align}
    \sigma'_z = \sigma_z + p = -0.12218 \, |z|^{1.0766} + 0.1|z| \label{eq:effectiveStress2Dtestcase}
\end{align}
\end{linenomath}
is the vertical effective stress, consisting of vertical total stress $\sigma_z$ and hydrostatic pressure $p$ (both in units of [$\text{bar}$]).   Young's modulus is expressed as
\begin{linenomath}
\begin{align}
    E = \frac{(1-2\nu)(1+\nu)}{(1-\nu)c_M} \label{eq:Youngs2Dtestcase}
\end{align}
\end{linenomath}
with a Poisson ratio $\nu$ set to 0.3 everywhere.   In addition, a constant value is added to Young's modulus in each layer, specifically the mean Young's modulus in that layer multiplied by a layer--dependent coefficient. This is done to emulate discontinuities in material properties between layers as often encountered in real subsurface systems.   The resulting distribution spans 3 orders of magnitude over the domain and is depicted in Figure~\ref{fig:cross_section_young_2D}.

\begin{figure} [htbp]
    \begin{subfigure}{0.48\textwidth}
        \includegraphics[width=\linewidth]{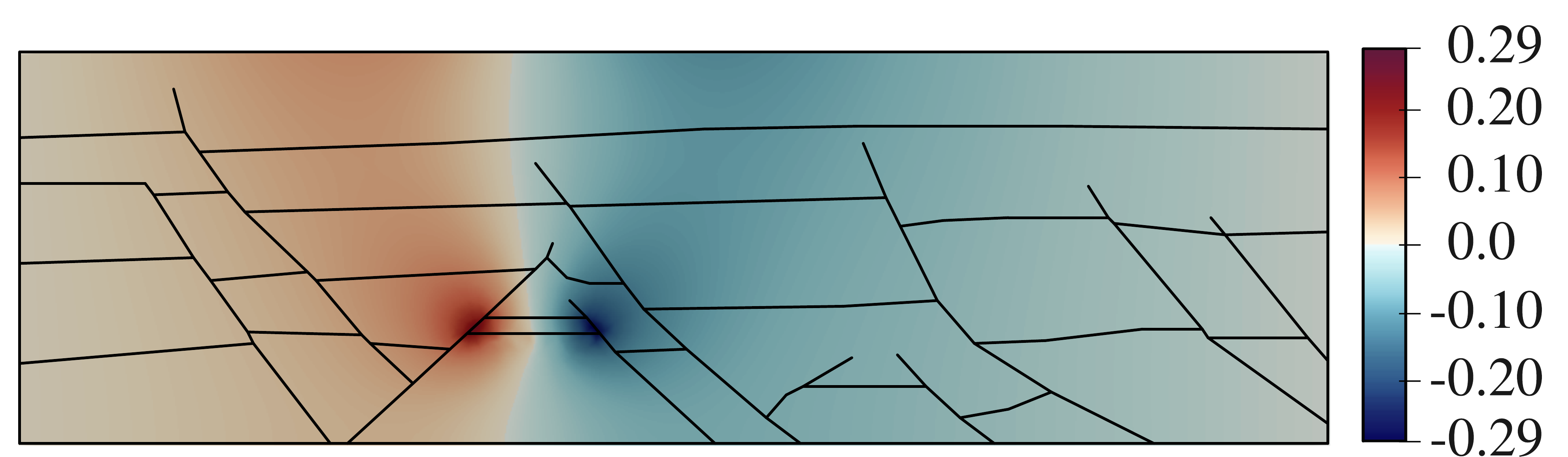}
	    \caption{Reference solution $x$-displacement [m]}
        \label{fig:cross_section_soln_ref_x}
    \end{subfigure}
    \hfill
    \begin{subfigure}{0.48\textwidth}
        \includegraphics[width=\linewidth]{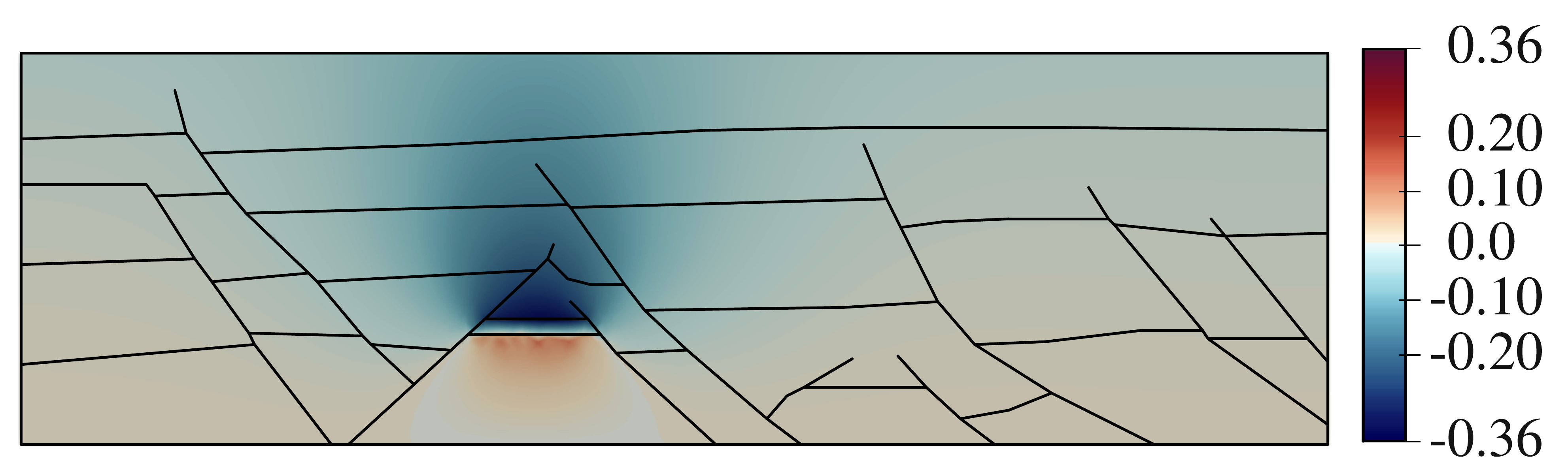}
	    \caption{Reference solution $z$-displacement [m]}
        \label{fig:cross_section_soln_ref_y}
    \end{subfigure} \\
    \begin{subfigure}{0.48\textwidth}
        \includegraphics[width=\linewidth]{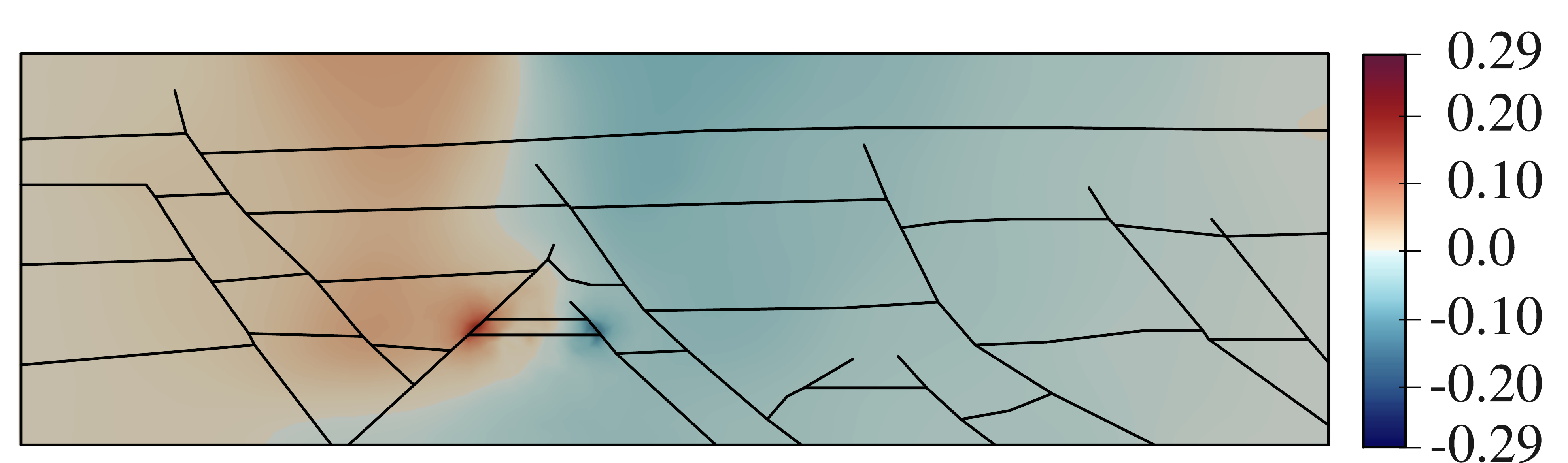}
	    \caption{Initial multiscale solution $x$-displacement [m]}
        \label{fig:cross_section_soln_ms_x}
    \end{subfigure}
    \hfill
    \begin{subfigure}{0.48\textwidth}
        \includegraphics[width=\linewidth]{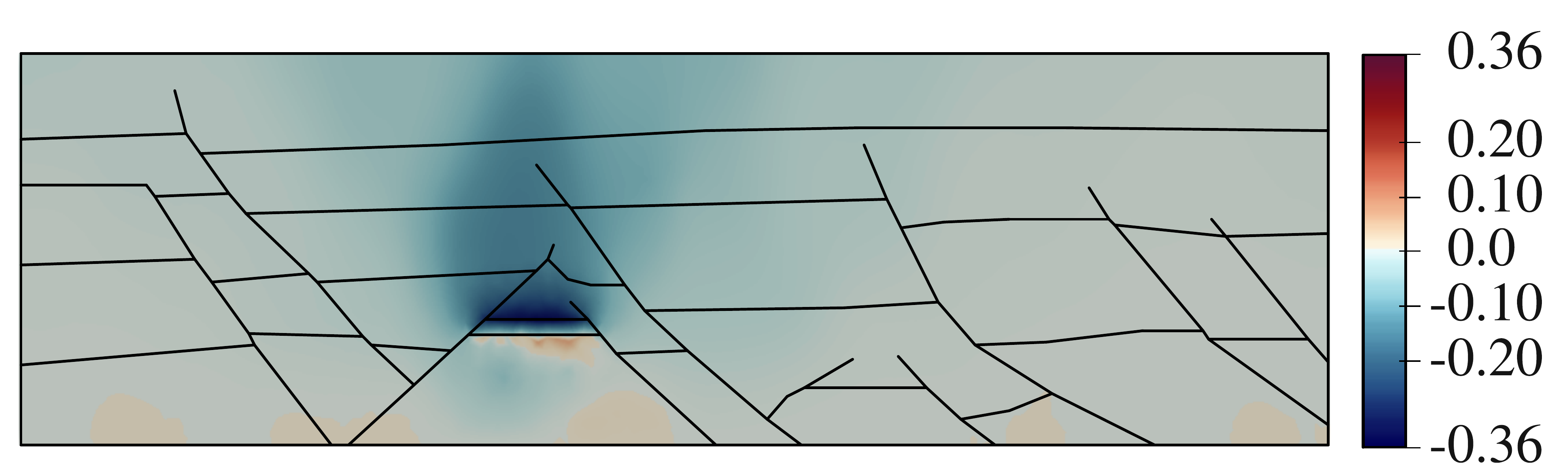}
	    \caption{Initial multiscale solution $z$-displacement [m]}
        \label{fig:cross_section_soln_ms_y}
    \end{subfigure} \\
    \begin{subfigure}{0.48\textwidth}
        \includegraphics[width=\linewidth]{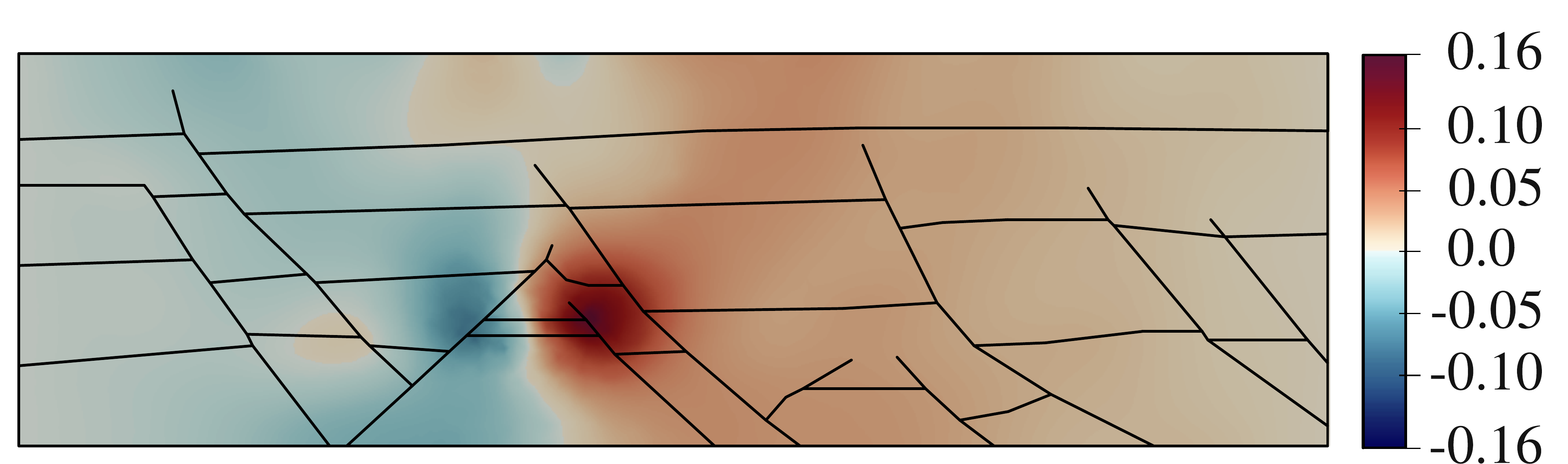}
	    \caption{Error in initial mutiscale solution $x$-displacement [m]}
        \label{fig:cross_section_error_abs_x}
    \end{subfigure}
    \hfill
    \begin{subfigure}{0.48\textwidth}
        \includegraphics[width=\linewidth]{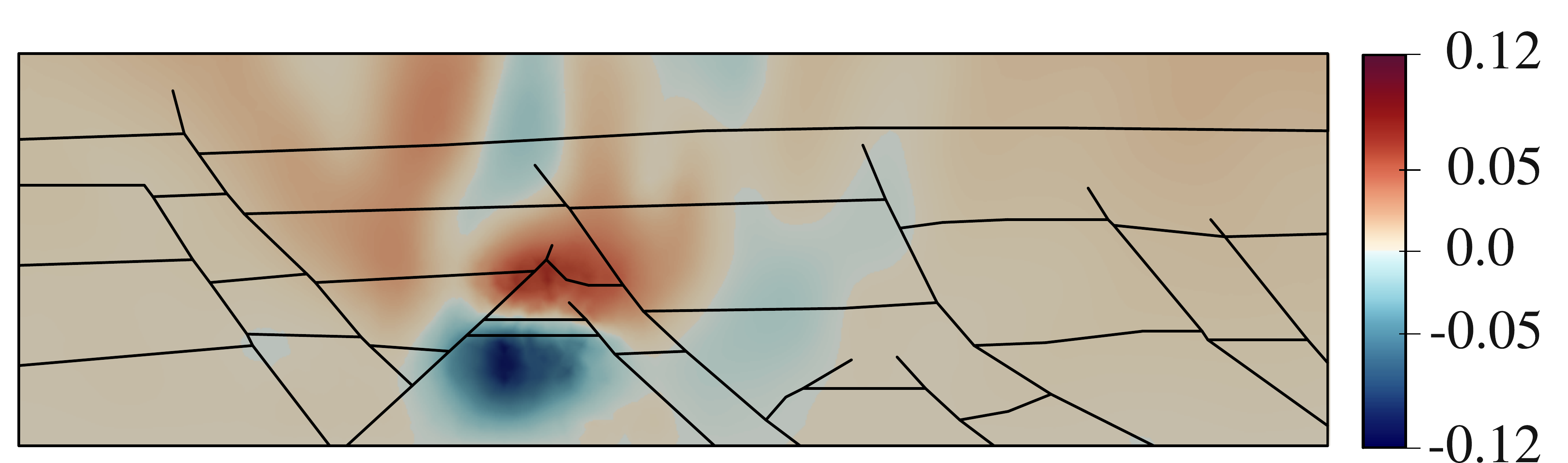}
	    \caption{Error in initial mutiscale solution $z$-displacement [m]}
        \label{fig:cross_section_error_abs_y}
    \end{subfigure}
    \caption{Geological 2D cross-section test case: reference solution (a,b) together with the initial, uniterated multiscale solution (c,d) and the corresponding error (e,f).}
    \label{fig:cross_section_soln}
\end{figure}

The domain is gridded with an unstructured triangular mesh that conforms to the layers and faults (see Figure~\ref{fig:cross_section_mesh}), with the faults themselves considered inactive (no slip between fault surfaces).   Seven different resolutions of the mesh are considered, ranging between 11,879 and 745,900 elements --- see Table~\ref{tab:cross_section_dims} for detailed information on mesh resolution and corresponding problem sizes.   The domain is subject to roller boundary conditions on three sides, while the ground surface is traction-free.   The deformation process is driven by a constant pressure drawdown of $\Delta p = 20$ bar prescribed in a small reservoir zone inside the domain (shown in Figure~\ref{fig:cross_section_sketch}), that acts as an external distributed force. Figure~\ref{fig:cross_section_soln} displays the the reference fine-scale displacement solution (computed on a grid corresponding to resolution level 0 in Table~\ref{tab:cross_section_dims}) in comparison with an approximate solution obtained from a single application of the two-stage multiscale preconditioner. Note that the smoothing stage is used to capture the effect of the forcing term prescribed in the interior of the domain via reservoir pressure drop.

To evaluate the algorithmic scalability of the MsRSB method for mechanics, the problem is solved with a Krylov solver (namely PCG) to a relative tolerance of $10^{-8}$.   Three symmetric two-stage preconditioning operators are constructed using different choices of pre- and post-smoother: 2 sweeps of $l_1$-Jacobi, symmetric Gauss-Seidel, and no-fill incomplete Cholesky factorization.   The coarse grid for the multiscale solver is generated by agglomerating fine-scale cells based on face connectivity using METIS \cite{metis} graph partitioning software.   For each mesh resolution, the ratio of fine--to--coarse elements and nodes is kept approximately the same, which results in the size of the coarse problem growing with mesh resolution.   As a comparative baseline, a smoothing-only preconditioner (i.e. not involving the global multiscale step) is also applied to the problem for each choice of smoother.   Krylov iteration counts are recorded to evaluate performance.

Table~\ref{tab:cross_section_conv} summarizes the findings.   The multiscale solver convergence remains well bounded for all mesh resolutions and only exhibits very mild mesh dependence, while the baseline approach does not scale well, in some cases failing to achieve convergence within 1000 CG iterations.   This example demonstrates a clear benefit of using a multiscale approach for subsurface mechanical problems compared to relying on incomplete factorizations only. It also emphasizes the method's excellent algorithmic scalability and mesh independence on unstructured grids.

\begin{table*}
    \centering
    \caption{Geological 2D cross--section test case: refinement levels and corresponding fine-- and coarse--scale problem sizes.}
    \label{tab:cross_section_dims}
    \begin{small}
    \begin{tabular}{rrrrrrrrrrcc}
        \hline\noalign{\smallskip}
        \multirow{2}{*}{$\ell$} & & \multicolumn{3}{c}{fine grid} & & \multicolumn{3}{c}{coarse grid} & & \multicolumn{2}{c}{coarsening ratio} \\
        \noalign{\smallskip}\cline{3-5} \cline{7-9} \cline{11-12}\noalign{\smallskip}
        & & \# cell & \# node & \# dof & & \# cell & \# node & \# dof & & cell & dof \\
        \hline\noalign{\smallskip}
        0 & &  11,879 &   6,119 &  12,238 & &  12 &    26 &    52 & & 989.9 & 235.4 \\
        1 & &  23,390 &  11,947 &  23,894 & &  25 &    51 &   102 & & 935.6 & 234.3 \\
        2 & &  46,932 &  23,817 &  47,634 & &  50 &   102 &   204 & & 938.6 & 233.5 \\
        3 & &  93,129 &  47,085 &  94,170 & & 100 &   202 &   404 & & 931.3 & 233.1 \\
        4 & & 186,940 &  94,165 & 188,330 & & 200 &   408 &   816 & & 934.7 & 230.8 \\
        5 & & 372,360 & 187,210 & 374,420 & & 400 &   803 & 1,606 & & 930.9 & 233.1 \\
        6 & & 745,900 & 374,350 & 748,690 & & 800 & 1,603 & 3,206 & & 932.4 & 233.5 \\
        \hline\noalign{\smallskip}
    \end{tabular}
    \end{small}
\end{table*}

\begin{table*}
    \centering
    \caption{Geological 2D cross--section test case: PCG iteration counts. ''MsRSB'' denotes the two-stage multiscale preconditioner with pre- and post-smoothing. ''no MS'' columns reports results for the baseline approach (only pre- and post-smoother, without the multiscale step). Dashes are shown for cases where the solver failed to converge in 1000 iterations. }
    \label{tab:cross_section_conv}
    \begin{small}
    \begin{tabular}{cccccccccc}
        \hline\noalign{\smallskip}
        \multirow{2}{*}{$\ell$} & & \multicolumn{2}{c}{$l_1$-Jacobi ($\times2$)} & & \multicolumn{2}{c}{Sym. Gauss-Seidel} & & \multicolumn{2}{c}{IC(0)} \\
        \noalign{\smallskip}\cline{3-4} \cline{6-7} \cline{9-10}\noalign{\smallskip}
        & & MsRSB & no MS & & MsRSB & no MS & & MsRSB & no MS \\
        \hline\noalign{\smallskip}
        0 & & \textbf{117} & 259 & & \textbf{66} & 145 & & \textbf{47} & 97  \\
        1 & & \textbf{113} & 369 & & \textbf{63} & 204 & & \textbf{43} & 137 \\
        2 & & \textbf{122} & 517 & & \textbf{68} & 290 & & \textbf{47} & 194 \\
        3 & & \textbf{126} & 727 & & \textbf{71} & 407 & & \textbf{48} & 268 \\
        4 & & \textbf{129} & --- & & \textbf{72} & 570 & & \textbf{49} & 385 \\
        5 & & \textbf{133} & --- & & \textbf{75} & 816 & & \textbf{51} & 544 \\
        6 & & \textbf{134} & --- & & \textbf{75} & --- & & \textbf{50} & 769 \\
        \hline\noalign{\smallskip}
    \end{tabular}
    \end{small}
\end{table*}

\subsection{MsRSB for Geomechanics: A 3D Test Case}
\begin{figure}[htbp]
    \begin{subfigure}[t]{0.3\textwidth}
        \centering
        \includegraphics[align=t,width=\textwidth]{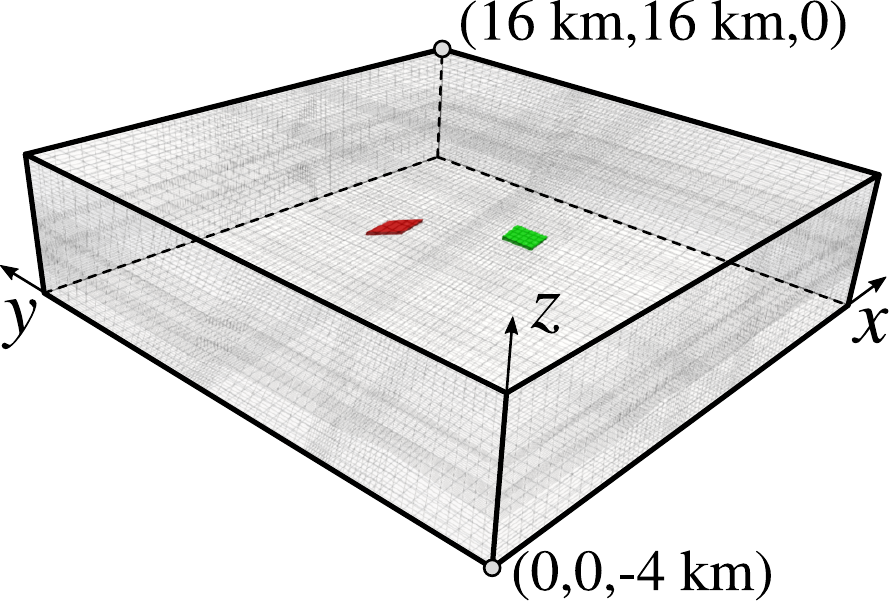}
        \caption{Mesh and reservoir zones.}
        \label{fig:3D_skew_mesh}
    \end{subfigure}
    \hfill
    \begin{subfigure}[t]{0.3\textwidth}
        \centering
        \includegraphics[align=t,width=\textwidth]{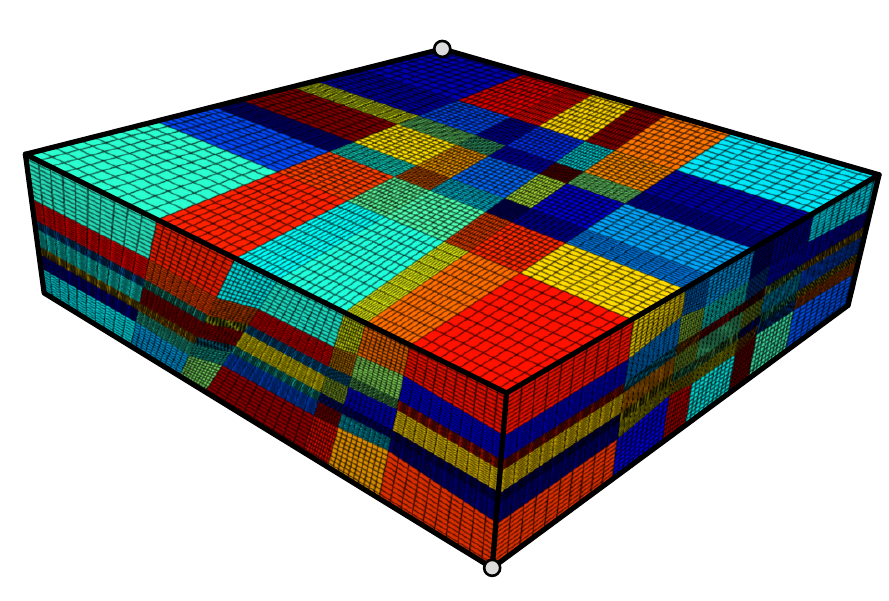}
        \caption{$7\times7\times7$ coarse-scale grid.}
        \label{fig:3D_skew_7x7x7}
    \end{subfigure}
    \hfill
    \begin{subfigure}[t]{0.3\textwidth}
        \centering
        \includegraphics[align=t,width=\textwidth]{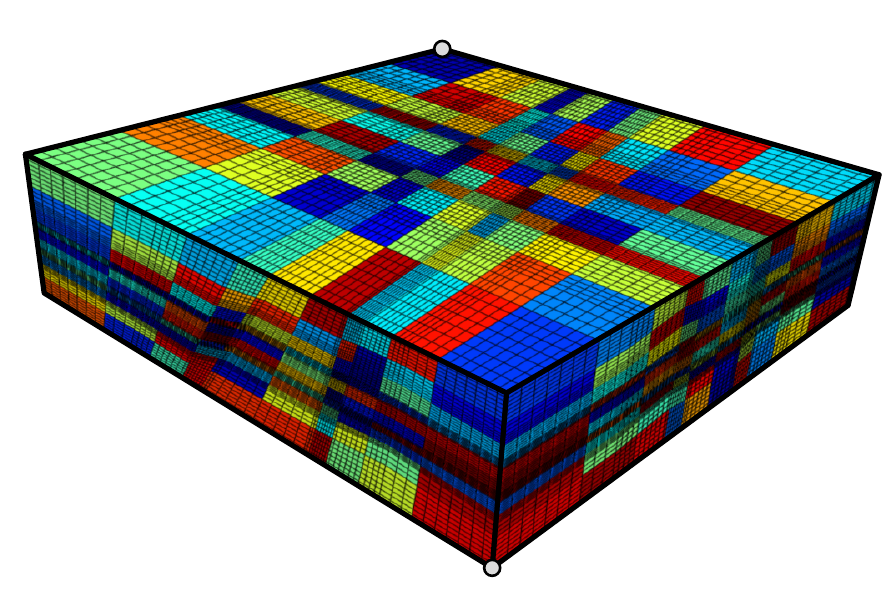}
        \caption{$10\times10\times10$ coarse-scale grid.}
        \label{fig:3D_skew_10x10x10}
    \end{subfigure}
    \caption{3D skewed mesh test case.}
    \label{fig:3D_structured_young}
\end{figure}

A third test is performed on a 3D poromechanical domain representing a $16\times16\times4$ km subsurface formation. The vertical distribution of Young's modulus prescribed by Eq. \ref{eq:compressibility2Dtestcase}--\ref{eq:Youngs2Dtestcase} is applied without discontinuities, and Poisson ratio is again set to 0.3. Also similar to the previous case, boundary conditions are imposed to be rollers (zero normal displacement) on all sides except for the traction-free top surface. Forcing is prescribed through pressure drawdown of $\Delta p_1 = 15$ bar and $\Delta p_2 = 22$ bar, respectively, in two reservoirs located around the center of the domain highlighted in Figure~\ref{fig:3D_skew_mesh}.

\begin{figure} [htbp]
\centering

\hfill
\begin{subfigure}{0.48\textwidth}
  \begin{tikzpicture}
    \begin{semilogyaxis}[ width  = \linewidth,
                          height = 1.\linewidth,
                          grid = major,
                          major grid style={very thin,draw=gray!50},
                          xmin=0,xmax=300,
                          ymin=1e-9,ymax=1,
                          xlabel={Iteration number},
                          ylabel={Relative residual norm},
                          xtick={0,50,...,300},
                          ytick={1e-9,1e-8,1e-7,1e-6,1e-5,1e-4,1e-3,1e-2,1e-1,1},
                          ylabel near ticks,
                          xlabel near ticks,
                          tick label style={font=\footnotesize},
                          label style={font=\footnotesize},
                          legend style={font=\footnotesize},
                          legend cell align={left},
                          legend columns=2,
                          legend pos=north east]

      \addplot [dashed, mygreen, line width=1pt]
               table[x index=0, y index=1, col sep=comma]
               {reservoir3dGrid_cart_7x7x7_CG_smoother_JACx2.csv};
      \addlegendentry{Jacobi};      
      \addplot [solid, mygreen, line width=1pt]
               table[x index=0, y index=1, col sep=comma]
               {reservoir3dGrid_cart_7x7x7_CG_multiscale_JACx2.csv}; 
      \addlegendentry{MsRSB + Jacobi};

      \addplot [dashed, blue, line width=0.8pt]
               table[x index=0, y index=1, col sep=comma]
               {reservoir3dGrid_cart_7x7x7_CG_smoother_SGSx1.csv};
      \addlegendentry{SGS};      
      \addplot [solid, blue, line width=0.8pt]
               table[x index=0, y index=1, col sep=comma]
               {reservoir3dGrid_cart_7x7x7_CG_multiscale_SGSx1.csv}; 
      \addlegendentry{MsRSB + SGS};
      
      \addplot [dashed, red, line width=0.8pt]
               table[x index=0, y index=1, col sep=comma]
               {reservoir3dGrid_cart_7x7x7_CG_smoother_IC0x1.csv};
      \addlegendentry{IC(0)};  
      \addplot [solid, red, line width=0.8pt]
               table[x index=0, y index=1, col sep=comma]
               {reservoir3dGrid_cart_7x7x7_CG_multiscale_IC0x1.csv}; 
      \addlegendentry{MsRSB + IC(0)};

    \end{semilogyaxis}
  \end{tikzpicture}
  \caption{Cartesian structured grid}
\end{subfigure}
\hfill
\begin{subfigure}{0.48\textwidth}
  \begin{tikzpicture}
    \begin{semilogyaxis}[ width  = \linewidth,
                          height = 1.\linewidth,
                          grid = major,
                          major grid style={very thin,draw=gray!50},
                          xmin=0,xmax=300,
                          ymin=1e-9,ymax=1,
                          xlabel={Iteration number},
                          ylabel={Relative residual norm},
                          xtick={0,50,...,300},
                          ytick={1e-9,1e-8,1e-7,1e-6,1e-5,1e-4,1e-3,1e-2,1e-1,1},
                          ylabel near ticks,
                          xlabel near ticks,
                          tick label style={font=\footnotesize},
                          label style={font=\footnotesize},
                          legend style={font=\footnotesize},
                          legend cell align={left},
                          legend columns=2,
                          legend pos=north east]

      \addplot [dashed, mygreen, line width=1pt]
               table[x index=0, y index=1, col sep=comma]
               {reservoir3dGrid_skew_7x7x7_CG_smoother_JACx2.csv};
      \addlegendentry{Jacobi};      
      \addplot [solid, mygreen, line width=1pt]
               table[x index=0, y index=1, col sep=comma]
               {reservoir3dGrid_skew_7x7x7_CG_multiscale_JACx2.csv}; 
      \addlegendentry{MsRSB + Jacobi};

      \addplot [dashed, blue, line width=0.8pt]
               table[x index=0, y index=1, col sep=comma]
               {reservoir3dGrid_skew_7x7x7_CG_smoother_SGSx1.csv};
      \addlegendentry{SGS};      
      \addplot [solid, blue, line width=0.8pt]
               table[x index=0, y index=1, col sep=comma]
               {reservoir3dGrid_skew_7x7x7_CG_multiscale_SGSx1.csv}; 
      \addlegendentry{MsRSB + SGS};
      
    \end{semilogyaxis}
  \end{tikzpicture}
  \caption{Skewed structured grid}
\end{subfigure}
\hfill\null

\caption{ CG-accelerated iterative performance for the Cartesian 3D mesh (a) and for the skewed structured 3D mesh (b).}
\label{fig:3D_structured_results}
\end{figure}
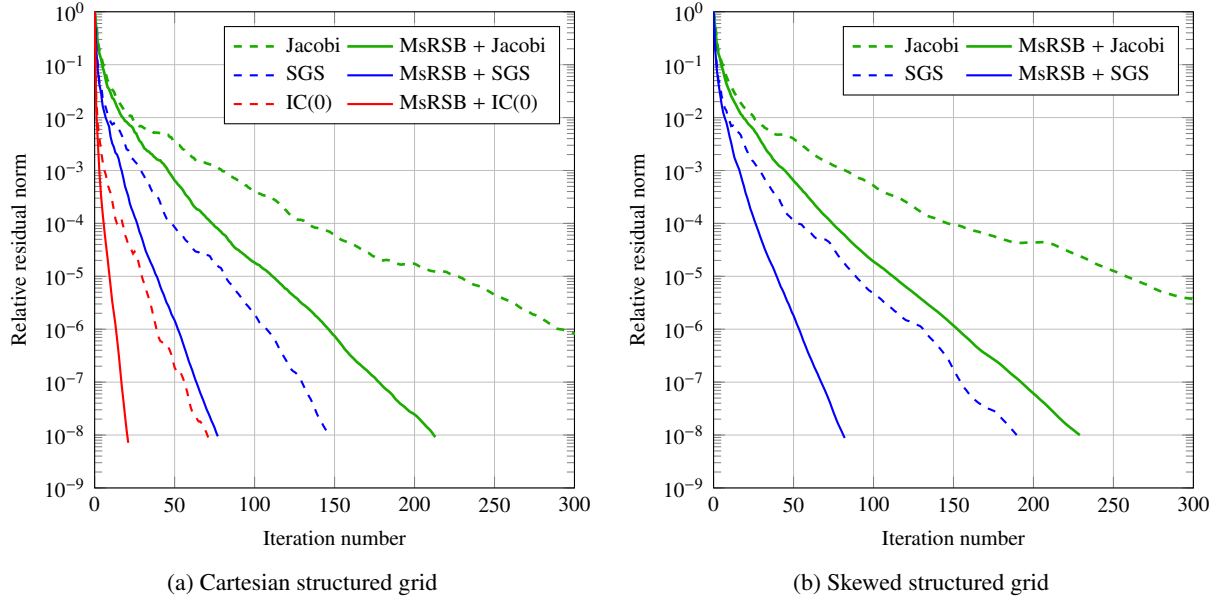

The domain is initially gridded with a $70\times70\times70$ structured Cartesian grid, labeled as \texttt{cart}. It is then coarsened with a sequence of progressively finer coarsening ratios leading to grids containing between 5 and 14 coarse cells in each dimension (examples are shown in Figures~\ref{fig:3D_skew_7x7x7}--\ref{fig:3D_skew_10x10x10}). For example a coarsening factor of 10 in each dimension, results in a $7\times7\times7$ coarse-scale grid. In addition, nodes of the grid are shifted, resulting in a skewed grid, shown in Figure~\ref{fig:3D_skew_mesh}, labeled \texttt{skew}.   For both grids, the multiscale preconditioner was constructed using the same choices of pre- and post-smoothers as in the previous example, i.e. $l_1$-Jacobi, symmetric Gauss-Seidel and no-fill incomplete Cholesky (the latter not being used with the skewed mesh due to numerical breakdowns in factorization which are not related to the multiscale method).  Table~\ref{tab:3D_structured_results} reports the observed iteration counts using CG as the chosen Krylov method and Figure~\ref{fig:3D_structured_results} compares convergence histories obtained using both multiscale (with a $10\times10\times10$ coarse element size) and smoother-only preconditioned CG solvers. 
The results are in line with the lower dimensional test cases. Acceptable iteration counts and good scalability is observed.

\begin{table*}
    \centering
    \caption{3D domain test case: coarse grid dimensions and CG iteration counts with different smoother options. }
    \label{tab:3D_structured_results}
    \begin{small}
    \begin{tabular}{ccccccccccc}
        \hline\noalign{\smallskip}
        \multirow{2}{*}{\# coarse cells} & \multirow{2}{*}{coarse cell size} & & \multicolumn{2}{c}{$l_1$-Jacobi $\times2$} & & \multicolumn{2}{c}{Gauss-Seidel} & & \multicolumn{2}{c}{IC(0)} \\
        \noalign{\smallskip}\cline{4-5} \cline{7-8} \cline{10-11}\noalign{\smallskip}
        & & & \texttt{cart} & \texttt{skew} & & \texttt{cart} & \texttt{skew} & & \texttt{cart} & \texttt{skew} \\
        \hline\noalign{\smallskip}
        $5\times5\times5$    & $14\times14\times14$ & & 257 & 269 & & 93 & 96 & & 24 & --- \\
        $7\times7\times7$    & $10\times10\times10$ & & 213 & 229 & & 77 & 82 & & 21 & --- \\
        $10\times10\times10$ & $7\times7\times7$    & & 178 & 203 & & 67 & 73 & & 17 & --- \\
        $14\times14\times14$ & $5\times5\times5$    & & 149 & 169 & & 58 & 63 & & 13 & --- \\
        \hline\noalign{\smallskip}
    \end{tabular}
    \end{small}
\end{table*}

\section{Conclusion}
A novel preconditioner based on Multiscale Restricted Smoothed Basis(MsRSB) functions is presented.
The essence of the enhanced MsRSB approach 
consists of enforcing M-matrix properties on the fine-scale linear system based on a filtering technique.
Using the resulting approximate linear system, basis functions can robustly be constructed using the original iterative MsRSB strategy. 

The method is demonstrated for the single phase flow problem and the linear elastic geomechanics problem. Enhanced MsRSB is validated to be effective through various test cases including heterogeneous, unstructured 2-dimensional and 3-dimensional problems . The proposed preconditioner show similar iteration counts compared to existing multiscale methods while enabling increased flexibility, easier implementation and sparser systems. Studies on CPU times and computational efficiency are subject of further work.

Finally, as shown for the porous media problem, the proposed method allows for the application of multiscale methods to multipoint stencils. Noting that multiscale operators are themselves inherently multipoint, the adapted MsRSB enables multilevel multiscale. Such a development has not yet been achieved in literature and is the topic of current research.

\section*{Acknowledgements} 
SB is supported by a named Stanford Graduate Fellowship in Science and Engineering (SGF). 
Funding for SK and NC was provided by Total S.A. through the FC-MAELSTROM Project.  
OM is funded by VISTA, which is a basic research program funded by Equinor and conducted in close collaboration with The Norwegian Academy of Science and Letters.
Portions of this work were performed under the auspices of the U.S. Department of Energy by Lawrence Livermore National Laboratory under Contract DE-AC52-07-NA27344.
The authors also thank Dr. Igor Shovkun and Dr. Jacques Franc for their feedback and Prof. Hamdi Tchelepi for many helpful discussions.

\appendix 

\section{Model problems: Governing equations} 
\label{app:models_strong_form}
Let $\Omega \subset \mathbb{R}^{n_{sd}}$ and $\Gamma$ denote a domain occupied by a heterogeneous porous medium and its boundary, respectively, with $\tensorOne{n}_{\Gamma}$ the unit outward normal vector to $\Gamma$, $\tensorOne{x}$ the position vector in $\mathbb{R}^{n_{sd}}$, and $n_{sd}$ (= 2 or 3) the spatial dimension of the problem.

\subsection{Incompressible single-phase flow} 
\label{app:model_flow}
Let $p$ denote the pore pressure.
For the application of the boundary conditions, let $\Gamma$ be decomposed as $\Gamma = \overline{\Gamma^D \cup \Gamma^N}$, where $\Gamma^D \cap \Gamma^N = \emptyset$.
The strong form of the incompressible single-phase flow boundary value problem (BVP) may be formally stated as follows: given $q : \Omega \rightarrow \mathbb{R}$, $g_D : \Gamma^D \rightarrow \mathbb{R}$, and $g_N : \Gamma^N \rightarrow \mathbb{R}$, find $p : \overline{\Omega} \rightarrow \mathbb{R}$ such that
\begin{linenomath}
\begin{subequations}
\begin{align}
	\text{div} \; \tensorOne{w}(p) &= q,
	&& \text{in} \; \Omega
	&& \mbox{(pressure equation)},
	\label{eq:massBalanceS} \\
  p &= g_{D},
  && \text{on} \; \Gamma^D
  && \mbox{(prescribed boundary pressure)},
  \label{eq:massBalanceS_DIR}\\
  \tensorOne{w}(p) \cdot \tensorOne{n}_{\Gamma} &= g_N,
  && \text{on} \; \Gamma^N
  && \mbox{(prescribed boundary flux)},
  \label{eq:massBalanceS_NEU}\\
	\tensorOne{w}(p) &= - \tensorTwo{\Lambda} \cdot \text{grad} \; p,
	&& \text{in} \; \Omega
	&& \mbox{(Darcy's law)}.
	\label{eq:Darcy}
\end{align}
\label{eq:model_flow_BVP}\null
\end{subequations}
\end{linenomath}
Here, $\tensorOne{w}(p)$ is the Darcy velocity, $\tensorTwo{\Lambda} = (\tensorTwo{\kappa} / \mu)$ is the the rank-two tensor characterizing the diffusion properties of the medium, with $\tensorTwo{\kappa}$ the intrinsic symmetric positive definite permeability tensor and $\mu$ the fluid viscosity, which is assumed constant, and $q$ denotes a volumetric source term.

\subsection{Linear elastostatics} 
\label{app:model_elastostatics}

Let $\tensorOne{d} = \{ d_{\ell} \}_{\ell = 1}^{n_{sd}}$ be the displacement-vector, where $d_{\ell} = ( \tensorOne{e}_{\ell} \cdot \tensorOne{d} )$ are the displacement components with respect to the Euclidean basis $\{ \tensorOne{e} \}_{\ell = 1}^{n_{sd}} $ in $\mathbb{R}^{n_{sd}}$.
Let us consider $n_{sd}$ non-overlapping partitions of the domain boundary into two segments associated with Dirichlet, $\Gamma_{\ell}^D$, and Neumann boundary conditions, $\Gamma_{\ell}^N$, respectively, such that $\Gamma = \overline{\Gamma_{\ell}^D \cup \Gamma_{\ell}^N}$, with $\Gamma_{\ell}^D \cap \Gamma_{\ell}^N = \emptyset$, $\ell \in \{ 1, \ldots, n_{sd} \} $.
The strong form of the linear elastostatic BVP reads as: given $\tensorOne{b} : \Omega \rightarrow \mathbb{R}^3$, $g_{D,\ell} : \Gamma_{\ell}^D \rightarrow \mathbb{R}$, and $g_{N,\ell} : \Gamma_{\ell}^N \rightarrow \mathbb{R}$, find $\tensorOne{d} : \overline{\Omega} \rightarrow \mathbb{R}^3$ such that
\begin{linenomath}
\begin{subequations}
\begin{align}
	- \text{div} \; \tensorTwo{\sigma}(\tensorOne{d}) &= \tensorOne{b},
	&& \text{in} \; \Omega
	&& \mbox{(equilibrium equations)},
	\label{momentumBalanceS} \\
  d_{\ell} &= g_{D,\ell},
  && \text{on} \; \Gamma_{\ell}^D
  && \mbox{(prescribed boundary displacements)},
  \label{momentumBalanceS_DIR}\\
  \tensorTwo{\sigma}(\tensorOne{d}) : \left( \tensorOne{e}_{\ell} \otimes \tensorOne{n}_{\Gamma} \right) &= g_{N,\ell},
  && \text{on} \; \Gamma_{\ell}^N
  && \mbox{(prescribed boundary tractions)},
  \label{momentumBalanceS_NEU}\\
	\tensorTwo{\sigma}(\tensorOne{d}) &= \tensorFour{C} : \text{sym} (\text{grad} \; \tensorOne{d}),
	&& \text{in} \; \Omega
	&& \mbox{(generalized Hooke's law)},
	\label{stressStrain}  
\end{align}
\label{eq:elasticity_global}\null
\end{subequations}
\end{linenomath}
with $\ell \in \{ 1, \ldots, n_{sd} \}$.. Here, $\tensorTwo{\sigma}(\tensorOne{d})$ is the rank-2 stress-tensor, respectively, $\tensorFour{C}_{dr}$ is the rank-4 elasticity tensor, and $\tensorOne{b}$ is a body force.
In this work we will focus on isotropic linear elastic materials, hence only two independent elastic coefficients are required for the definition of $\tensorFour{C}$, namely
\begin{align}
  \tensorFour{C} &= \lambda (\tensorTwo{I} \otimes \tensorTwo{I}) + 2G \tensorFour{I},
  \label{eq:elasticity_tensor}
\end{align}
where $\tensorTwo{I}$ and $\tensorFour{I}$ are the second-order and fourth-order identity tensor, respectively, and $\lambda = \frac{E \nu}{(1+\nu)(1-2\nu)}$ and $G = \frac{E}{2(1+\nu)}$ are the Lam\'{e} parameters of the material, with $E$ the Young modulus and $\nu$ the Poisson ratio.
%
Note that the subscripts $x$, $y$, and $z$ are also used to denote a quantity associated with the spatial dimension $\ell$ equal to 1, 2, and 3, respectively.
    
\section{Model problems: Discrete formulation} 
\label{app:model_problems_discretization}
\subsection{Incompressible single-phase flow: Finite volume formulation} 
\label{app:model_flow_FV}

Given a partition $\mathcal{T}^h$ of the domain $\Omega$ consisting of non-overlapping conforming cells, a finite volume discretization of \eqref{eq:model_flow_BVP} consists of writing the pressure equation for each cell (control-volume) in $\mathcal{T}^h$ in integral form \cite{EymGalHer00}.
Let $\mathcal{F}^h$ be the set of interfaces in $\mathcal{T}^h$, namely edges ($n_{sd} = 2$), or faces ($n_{sd} = 3$).
Let $\funSpace{V}^h$ be the space of piecewise constant cell-wise fuctions associated with $\mathcal{T}^h$.
We consider a discrete approximation for the pressure field such that $p \approx u^h \in \funSpace{V}^h$.
Let $g_D^h$ denote the piecewise constant interface-wise interpolant of $g_D$ having support in $\mathcal{F}^{h,D} \subset \mathcal{F}^h$, namely the set of interfaces belonging to $\Gamma^D$.
Similarly, $\mathcal{F}^{h,N} \subset \mathcal{F}^h$ is the set of boundary interfaces lying on $\Gamma^N$.
Finally, let $\hat{w}^{\gamma}$ denote a conservative numerical flux approximating the volumetric flux through an interface $\gamma \in \mathcal{F}^h$, namely $\hat{w}^{\gamma} \approx \int_{\gamma} \tensorOne{w} \cdot \tensorOne{n}_{\gamma} \, \mathrm{d}\Gamma$, with $\tensorOne{n}_{\gamma}$ a unit normal vector defining a unique global orientation for $\gamma$.
Based on a suitable functional dependence on $u^h$ and $g_D^h$, for linear flux approximation schemes $\hat{w}^{\gamma}(u^h,g_D^h)$ can be split as a sum of two terms;

\begin{align}
  \hat{w}^{\gamma}(u^h,g_D^h) = \mathring{w}^{\gamma}(u^h) + \bar{w}^{\gamma}(g_D^h),
  \label{eq:num_flux_split}
\end{align}

\noindent
to highlight the contribution to the flux related to $u^h$ and $g_D^h$, respectively.
Clearly, $\bar{w}^{\gamma}(g_D^h)$ is nonzero only in the presence of non homogeneous pressure boundary conditions.
The first term to the right-hand side in \eqref{eq:num_flux_split} is expressed as a linear combination of pressure values from selected cells---e.g., the two cells sharing $\gamma$ in the TPFA method, or the cells sharing at least a vertex with $\gamma$ in the MPFA-O method \cite{EdwRog98,MPFA_Aavatsmark}---using constant transmissibility coefficients.
A similar linear expression is utilized for $\bar{w}^{\gamma}(g_D^h)$.
Nonlinear flux approximation schemes are not considered in this work.
For a review and details on recent developments on finite volume discretizations for anisotropic diffusion problems in heterogeneous media, we refer the reader to \cite{Dro14,TerMalTch17} and references therein.

The finite volume discretization provides an approximation $u^h$ to the weak pressure solution of\eqref{eq:model_flow_BVP} by solving a set of discrete balance equations that are equivalent to the following mesh-dependent variational problem: find $u^h \in \funSpace{V}^h$ such that

\begin{align}
  a^h( v, u) = F^h(v) \qquad \forall v^h \in \funSpace{V}^h,
  \label{eq:model_flow_weak_FV}
\end{align}

\noindent
where the discrete bilinear form $a^h : \funSpace{V}^h \times \funSpace{V}^h \rightarrow \mathbb{R}$ and the discrete linear form $F: \funSpace{V}^h \rightarrow \mathbb{R}$  are defined as

\begin{align}
  &a^h(v^h,u^h)
  =
  - \sum_{\gamma \in \mathcal{F}^h \setminus \mathcal{F}^{h,N}} \llbracket v^h \rrbracket_{\gamma} \mathring{w}^{\gamma}(u^h),
  \label{eq:model_flow_bilin_FV}\\
  &F^h(v^h)
  =
  \int_{\Omega} v^h q \, \mathrm{d}\Omega
  + \sum_{\gamma \in \mathcal{F}^{h,N}} \llbracket v^h \rrbracket_{\gamma} \int_{\gamma} \bar{q} \, \mathrm{d}\Gamma
  + \sum_{\gamma \in \mathcal{F}^h \setminus \mathcal{F}^{h,N}} \llbracket v^h \rrbracket_{\gamma} \bar{w}^{\gamma}( g_D^h).
  \label{eq:model_flow_functional_FV} 
\end{align}

\noindent
Here, the symbol $\llbracket \cdot \rrbracket_{\gamma}$ denotes the jump of a quantity across an interface $\gamma\in \mathcal{F}^h$.
For internal interfaces, $\llbracket v^h \rrbracket_{\gamma} = ( {v^h}_{\left| \tau_L \right.} - {v^h}_{\left| \tau_K \right.} )$, with ${v^h}_{\left| \tau_L \right.}$ and ${v^h}_{\left| \tau_K \right.}$ the restriction of $v^h$ on cells $\tau_K$ and $\tau_L$ sharing $\gamma$, with $\tensorOne{n}_{\gamma}$ pointing from $\tau_K$ to $\tau_L$.
For domain boundary interfaces, the jump expression simplifies to $\llbracket v^h \rrbracket_{\gamma} = - {v^h}_{\left| \tau_K \right.}$.

To obtain the matrix form of the FV discrete problem, we introduce the basis $\{ \chi_i \}_{i \in \mathcal{N}^h}$ for $\funSpace{V}^h$, with $\chi_i$ the characteristic function of the $i$th cell $\tau_i$ in $\mathcal{T}^h$ such that $\chi_i(\tensorOne{x}) = 1$, if $\tensorOne{x} \in \tau_i$, $\chi_i(\tensorOne{x}) = 0$, if $\tensorOne{x} \notin \tau_i$, and $\mathcal{N}^h = \{1, \ldots, n_{\tau} \}$ with $n_{\tau}$ the total number of cells.
Hence, the approximate pressure field is expressed as $p(\tensorOne{x}) \approx u^h(\tensorOne{x}) = \sum_{i \in \mathcal{N}^h} u_i \chi_i(\tensorOne{x})$, with $u_i$ the unknown cell pressure values.
Requiring that $u^h$ satisfy \eqref{eq:model_flow_functional_FV} for each function of the basis itself yields the system of discrete balance equations for the unknown coefficients vector $\Vec{u} = \{ u_i \}$
\begin{linenomath}
\begin{align}
  \Mat{A}_p \Vec{u} = \Vec{f}_p,
  \label{eq:model_flow_linsys_FV}
\end{align}
\end{linenomath}
with the system matrix $\Mat{A}_p$ and the right-hand side $\Vec{f}_p$ such that $[\Mat{A}_p]_{ij} = a( \chi_i, \chi_j)$ with $\{ i, j \} \in \mathcal{N}^h \times \mathcal{N}^h$ and $\{\Vec{f}_p\}_i = F( \chi_i)$, with $i \in \mathcal{N}^h$.
The properties of matrix $A_p$ depend on the flux approximation scheme chosen for $\hat{w}^{\gamma}$.
For example, a TPFA scheme produces a symmetric positive definite matrix whereas MPFA methods typically lead to a non-symmetric $A_p$ \cite{Dro14}.

\subsection{Linear elastostatics: Galerkin finite element formulation} 
\label{app:model_elasticity_FE}

The weak form of the linear elastostatics BVP is derived based on the classical displacement formulation~\cite{Hug00} by eliminating $\tensorTwo{\sigma}(\tensorOne{d})$ in \eqref{momentumBalanceS} using \eqref{stressStrain}.
Under appropriate regularity assumptions, \eqref{eq:model_flow_BVP} admits a unique solution $\tensorOne{d}$ that can be obtained by solving an equivalent variational problem \cite{Hug00}.
Let $\vecFunSpace{V} = \{ \tensorOne{v} \in [H^1(\Omega)]^{n_{sd}} : v_{\ell} = ( \tensorOne{e}_{\ell} \cdot \tensorOne{v} ) \in {V}_{\ell} \}$ denote the space of test functions, where ${V}_{\ell} = \{ v \in H^1(\Omega): v | _{\Gamma^D_{\ell}} = 0 \}$, $\ell \in \{1, \ldots, n_{sd} \}$.
Let us consider an extension of the Dirichlet boundary datum $\tilde{\tensorOne{g}}_D \in [H^1(\Omega)]^{n_{sd}}$ such that $( \tensorOne{e}_{\ell} \cdot \tilde{\tensorOne{g}}_D) =g_{D,\ell}$ on $\Gamma^D_{\ell}$, $\ell \in \{1, \ldots, n_{sd} \}$. 
By expressing the displacement as $\tensorOne{d} = \tilde{\tensorOne{g}}_D + \tensorOne{u}$, the weak form of \eqref{eq:elasticity_global} reads as: find $\tensorOne{u} \in \vecFunSpace{V}$ such that

\begin{align}
  a( \tensorOne{v}, \tensorOne{u}) = F(\tensorOne{v}) \qquad \forall \tensorOne{v} \in \vecFunSpace{V},
  \label{eq:model_elastostatics_weak}
\end{align}

\noindent
where the bilinear form $a : \vecFunSpace{V} \times \vecFunSpace{V} \rightarrow \mathbb{R}$ and the linear form $F: \vecFunSpace{V} \rightarrow \mathbb{R}$  are defined as

\begin{align}
  &a(\tensorOne{v}, \tensorOne{u})
  =
  \int_{\Omega} \text{sym} (\text{grad} \; \tensorOne{v}) : \tensorFour{C}_{dr} : \text{sym} (\text{grad} \; \tensorOne{u}) \, \mathrm{d} \Omega,
  \label{eq:model_elastostatics_bilin_FE}\\
  &F(\tensorOne{v})
  =
  \int_{\Omega} \tensorOne{v} \cdot \tensorOne{b} \, \mathrm{d}\Omega
  + \sum_{\ell = 1}^{n_{sd}} \int_{\Gamma^N_{\ell}} v_\ell g_{N,\ell} \, \mathrm{d}\Gamma
  - \int_{\Omega} \text{sym} (\text{grad} \; \tensorOne{v}) : \tensorFour{C}_{dr} : \text{sym} (\text{grad} \; \tilde{\tensorOne{g}}_D) \, \mathrm{d}\Omega
  \label{eq:model_elastostatics_functional_FE}.
\end{align}   

\noindent
Let $\vecFunSpace{X}^h$ be the finite element space of piecewise polynomial vector functions that are continuous in $\overline{\Omega}$ associated with a conforming triangulation $\mathcal{T}^h$ of $\Omega$.
Let $\{ \tensorOne{\eta}_i \}_{i \in \mathcal{N}^h}$ be the standard (vector) nodal basis for $\vecFunSpace{X}^h$, with $\mathcal{N}^h = \{1, \ldots, n_{sd}n_n \}$ and $n_n$ the number of node points in $\mathcal{T}^h$.
We define the finite dimensional counterpart of $\vecFunSpace{V}$ as $\vecFunSpace{V}^h = \vecFunSpace{X}^h \cap \vecFunSpace{V}$ and denote its basis as $\{ \tensorOne{\eta}_i \}_{i \in \mathcal{N}^h_u}$, with $\mathcal{N}^h_u \subset \mathcal{N}^h$ the set of indexes of basis functions of $\vecFunSpace{X}^h$ vanishing on $\Gamma^D_{\ell}$, $\ell \in \{1, \ldots, n_{sd} \}$.
The discrete approximation to the displacement field can then be expressed as

\begin{align}
  \tensorOne{d}(\tensorOne{x})
  \approx  \tilde{\tensorOne{g}}_D^h(\tensorOne{x})
  + \tensorOne{u}^h(\tensorOne{x})
  = \sum_{\ell = 1}^{n_{sd}}
  \sum_{j \in \mathcal{N}^h \setminus \mathcal{N}^h_u} g_{D,\ell}(\phi_{\ell j}) \phi_{\ell j}(\tensorOne{x}) \tensorOne{e}_{\ell}
  + \sum_{j \in \mathcal{N}^h_u} u_j \tensorOne{\eta}_j(\tensorOne{x}),
  \label{eq:approx_u_FE}  
\end{align}

\noindent
where $\tilde{\tensorOne{g}}_D^h \in \vecFunSpace{X}^h$ is the trivial discrete extension of the Dirichlet boundary datum such that $\tensorOne{e}_{\ell} \cdot \tilde{\tensorOne{g}}_D^h$ on $\Gamma^D_{\ell}$ is equal to the finite element interpolant of $g_{D,\ell}$, $\ell \in \{1, \ldots, n_{sd} \}$, $\phi_{\ell j} = (\tensorOne{e}_{\ell} \cdot \tensorOne{\eta}_j(\tensorOne{x}))$ , and $\tensorOne{u}^h \in \vecFunSpace{V}^h$ is an approximate solution to the corresponding homogeneous Dirichlet problem with $u_j$ the unknown nodal displacement degrees of freedom.  
Substituting $\tilde{\tensorOne{g}}_D$ by $\tilde{\tensorOne{g}}_D^h$  in \eqref{eq:model_elastostatics_functional_FE} and requiring that $\tensorOne{u}^h$ satisfy \eqref{eq:model_elastostatics_weak} for each basis function of $\vecFunSpace{V}^h$ yields the matrix form of the variational problem, namely the system of equations for the unknown coefficients vector $\Vec{u} = \{ u_j \}$

\begin{align}
  \Mat{A}_d \Vec{u} = \Vec{f}_d,
  \label{eq:model_elastostatics_linsys}
\end{align}

\noindent
with $\Mat{A}_d$ and $\Vec{f}_d$ the symmetric positive definite (SPD) stiffness matrix and force vector, respectively, such that $[\Mat{A}_d]_{ij} = a( \tensorOne{\eta}_i, \tensorOne{\eta}_j)$ with $\{ i, j \} \in \mathcal{N}^h_u \times \mathcal{N}^h_u$ and $\{\Vec{f}_d\}_i = F( \tensorOne{\eta}_i)$, with $i \in \mathcal{N}^h_u$.
%

\begin{rem}\label{rem:FEM_Dir_BC}
For efficiency reasons, the linear system \eqref{eq:model_elastostatics_linsys} is typically assembled ignoring the essential (Dirichlet) boundary conditions---i.e., $\{ \tensorOne{\eta}_i, \tensorOne{\eta}_j \}$ is the range over the bases for $\vecFunSpace{X}^h$.
This is also the strategy adopted in our implementation.
The Dirichlet conditions are introduced by using a so-called \textit{symmetric diagonalization} approach \cite{For_etal12}.
In this approach, rows and columns of the stiffness matrix associated with displacement degrees of freedom where such conditions apply are modified while preserving symmetry, with the right-hand-side updated accordingly.
\end{rem}

\begin{rem}\label{rem:SDC}
If displacement degrees of freedom are ordered based on each coordinate direction, $A_d$ possesses a $n_{sd} \times n_{sd}$ block structure, namely:

\begin{align}
	A_d &=
    \begin{bmatrix}
  	  A_{11} & \ldots & A_{1n_{sd}}\\
  	  \vdots & \ddots & \vdots\\
  	  A_{n_{sd}1} & \ldots & A_{n_{sd}n_{sd}}
	\end{bmatrix},
  \label{eq:blk_stiff_gen}
\end{align} 

\noindent
which reflects the full coupling between $\ell$-components of displacements, $\ell = \{1, \ldots, n_{sd} \}$.
Each diagonal block $A_{\ell\ell}$ in \eqref{eq:blk_stiff} corresponds to the finite element discretization of the anisotropic diffusion operator $-\text{div} (\tensorTwo{\Lambda} \cdot \text{grad} \; d_{\ell})$, with anisotropic diffusion tensor $\tensorTwo{\Lambda} = G \tensorTwo{I} + (G+\lambda)\  \tensorOne{e}_{\ell} \otimes \tensorOne{e}_{\ell} $, $\ell = \{1, \ldots, n_{sd} \}$.
\end{rem}

\bibliography{journalAbbreviationExtended,MsRSBgeom}


\end{document}